\documentclass[11pt,a4paper]{article}
\usepackage{amssymb,fullpage,amsmath,mathrsfs}
\usepackage{pgfplots}
\usepackage{amsfonts,amsbsy,graphicx}
\usepackage{authblk,lineno}
\usepackage{bm}

\newcommand{\eeref}[1]{(\ref{eqn:#1})}
\newcommand{\eelab}[1]{\label{eqn:#1}}
\newcommand{\ffref}[1]{\ref{fig:#1}}
\newcommand{\fflab}[1]{\label{fig:#1}}
\newcommand{\ssref}[1]{\ref{sec:#1}}
\newcommand{\sslab}[1]{\label{sec:#1}}

\def\beq{\begin{equation}}
\def\eeq{\end{equation}}

\DeclareMathOperator\erfc{erfc}
\DeclareMathOperator\erf{erf}

\def\XXint#1#2#3{{\setbox0=\hbox{$#1{#2#3}{\int}$}
\vcenter{\hbox{$#2#3$}}\kern-.5\wd0}}

\def\strutdepth{\dp\strutbox}
\def\nw#1{\strut\vadjust{\kern-\strutdepth\vtop to0pt{\vss\hbox to\hsize {\hskip\hsize\hskip5pt$\leftarrow$\hss\strut}}}{\em \textcolor{blue}{#1}}}

\begin{document}

\title{The Evolution of Travelling Waves in a KPP Reaction-Diffusion Model with cut-off Reaction Rate.  II. Evolution of Travelling Waves.}

\author{A. D. O. Tisbury, D. J. Needham and A. Tzella\thanks{Address for correspondence:   Prof. D. J. Needham and Dr A. Tzella, School of Mathematics, University of Birmingham; email:  a.tzella@bham.ac.uk}} 
\affil{School of Mathematics, University of Birmingham, Birmingham, B15 2TT, UK.}
\maketitle

\begin{abstract}
In Part II of this series of papers, we consider an initial-boundary value problem for the Kolmogorov--Petrovskii--Piscounov (KPP) type equation with a discontinuous cut-off in the reaction function at concentration $u=u_c$. For  fixed cut-off value $u_c \in (0,1)$, we apply the method of matched asymptotic coordinate expansions to obtain the complete   large-time asymptotic form of the solution which exhibits the formation of a permanent form travelling wave structure. In particular, this approach allows the correction to the wave speed and the rate of convergence of the solution onto the permanent form travelling wave to be determined via a detailed analysis of the asymptotic structures in small-time and, subsequently, in large-space. The asymptotic results are confirmed against numerical results obtained for the particular case of a cut-off Fisher reaction function.
\end{abstract}

\section{Introduction}
Travelling waves arise as the long-time solution
to  many reaction-diffusion models and are relevant to a broad range of applications in chemistry, biology, ecology, epidemiology and genetics  \cite{Fife1979,Murray2002}. The  most celebrated model where such waves emerge is the KPP or Fisher-KPP model
named after the pioneering work by Fisher \cite{Fisher1937}
and Kolmogorov, Petrovskii, Piscounov \cite{Kolmogorov_etal1937}. 
In one spatial coordinate ($x$) this model describes the temporal ($t$) evolution of the concentration of a chemical or biological substance  
 $u(x,t)$ as 
 \begin{subequations}\eelab{KPP}
 \begin{align}
 & u_t  = u_{xx} + f(u), \quad (x,t) \in \mathbb{R} \times \mathbb{R}^+,\eelab{KPPa}\\
 \intertext{subject to an   initial condition  }
 & u(x,0)=u_0(x), \quad x \in \mathbb{R} \eelab{KPPb} \\
\intertext{
and  boundary conditions}
  & u(x,t)  \to
	 \begin{cases}
		 1, & \text{as $x \to - \infty$}\\
	             0, & \text{as $x \to \infty$},
	\end{cases} \eelab{KPPc}
  \end{align}
 \end{subequations}
  with the limits being uniform for time $t\in[0,T]$ and any $T>0$.
Here, $u_0:\mathbb{R}\to\mathbb{R}$  is taken to be 
piecewise continuous, non-negative and non-increasing with 
$\lim_{x\to\infty}u_0(x)=0$ and $\lim_{x\to-\infty}u_0(x)=1$. 
The  function
  $f: \mathbb{R} \to \mathbb{R}$  is a normalised KPP-type reaction function which satisfies  $f \in C^1(\mathbb{R})$ 
with 
  \begin{subequations} \eelab{KPPreaction}
 	\beq
	  	f(0) = f(1)= 0,\quad f'(0)  = 1, \quad f'(1)  <0
	\eeq	
	 	and 
		\beq
	  	 0<f(u)\leq u   \quad  \text{for all}\quad    u \in (0,1),\quad
	 	\quad f(u) < 0 \quad  \text{for all}\quad    u \in (1,\infty).
	  	\eeq
	  \end{subequations}  
	  A prototypical example of such a KPP reaction function is the Fisher reaction function \cite{Fisher1937} given by 
	  \beq\eelab{Fisher}
	  	  f(u)=u(1-u).
	  \eeq
	  The initial-boundary value problem  \eeref{KPP} has a classical  and global solution $u: \mathbb{R}\times[0,\infty)\to\mathbb{R}$. In addition, 
	  on using the classical 
	  maximum principle and comparison theorem (see, for example, \cite{AronsonSerrin1967} and \cite{Fife1979}),
	  $0 < u(x,t) < 1$ and $u_x(x,t)<0$ for all $(x,t) \in\mathbb{R}\times \mathbb{R}^+$.
The conditions \eeref{KPPreaction} on $f$ imply also that 
the initial-boundary value problem \eeref{KPP}
admits a one-parameter family of permanent form travelling wave (PTW) solutions  
$u(x,t)=U_v(x-vt)$ that are    strictly monotone decreasing, with $U_v\geq 0$,
$U_v: \mathbb{R}:\mathbb{R}$ such that $U_v>0$ with  
 $\lim_{y\to-\infty}U_v(y)=1$ 
and $\lim_{y\to\infty}U_v(y)=0$. 
The parameterisation is through the propagation speed $v$, with a unique 
(up to translation) PTW for each $v$ where $v$ satisfies 
 $v\geq v_m=2$.

A central question is whether
a  PTW 
evolves  in the
solution to  \eeref{KPP}  
 at large times  
and if so what is its speed of propagation. 
It is well established \cite{AronsonWeinberger1975,FifeMcLeod1977,Kolmogorov_etal1937} that for Heaviside initial conditions:
\beq\eelab{step}
 u_0=\begin{cases}
		 1, & \text{for $x <0$}\\
	             0, & \text{for $x \geq 0$},
	\end{cases} 
\eeq
 the solution to  \eeref{KPP}  
 converges onto
the PTW  solution  with minimum propagation speed $v=v_m=2$
in the sense that there exists a function $s_m(t)$
such that as $t\to\infty$,
 $s_m(t)/t \to 2$ and 
\beq
u(z+s_m(t),t)\to U_2(z),
 \eeq
uniformly for $z\in\mathbb{R}$. 
A more detailed asymptotic description was 
 provided by McKean \cite{McKean1975,McKean1976} and Bramson \cite{Bramson1978,Bramson1983}
who, using a probabilistic approach, obtained 
that the rate of 
 convergence of  the solution to the initial-boundary value problem \eeref{KPP} to the 
PTW is  algebraically small  in $t$
as $t\to\infty$, specifically
$O(\dot s_m(t)-2)$, 
where 
\beq\eelab{BramsonFisher}
\dot s_m(t)=2-\frac{3}{2}t^{-1}+o(t^{-1})\quad\text{as $t\to\infty$}
\eeq
with the dot denoting differentiation with respect to $t$.
More recently, the same result has been established  
using a range of alternative approaches, based on  a point patching procedure \cite{BrunetDerrida1997,EbertvanSarloos2000}, 
the theory of matched asymptotic expansions \cite{BillinghamNeedham1992,LeachNeedham2003} and 
rigorous bounds  \cite{Hamel_etal2013}. 
All of these approaches  involve the solution to a linearized version  of  \eeref{KPP} that  describes the behaviour at the leading edge
of the front and is 
obtained by replacing $f(u)$ with $f'(0)u$. The common observation is that, with the appropriate boundary conditions, 
the linear version  of  \eeref{KPP} mainly determines 
the large-$t$ structure of the solution to \eeref{KPP}. 

A linearized approach 
is not available to apply  
in the case of the 
cut-off KPP model that Brunet and Derrida \cite{BrunetDerrida1997}
proposed and considered and was the focus of a companion paper \cite{Tisbury_etal} (hereafter referred to as Part I). 
In this model, the cut-off value $u_c\in(0,1)$   is introduced by replacing 
$f(u)$ in the initial-boundary value problem \eeref{KPP} with $f_c(u)$ where 
\beq \eelab{BDreaction}
f_c(u)=
\begin{cases}
 			f(u), & \text{$u \in (u_c, \infty)$}\\
             0, & \text{$u \in (- \infty, u_c]$} \end{cases}
\eeq 
and $f(u)$ continues to satisfy the KPP conditions \eeref{KPPreaction}. 
The discontinuity in $f_c(u)$ at $u=u_c$  suggests that 
the  corresponding initial-boundary value problem  is expressed  as a moving boundary  problem 
with the location of the moving boundary given by 
 $s(t)$  where $s(t)$ satisfies $u(s(t),t)=u_c$ for $t>0$  
 (see Part I). For Heaviside initial conditions \eeref{step},
this boundary separates the  domain $D^L$ where $u>u_c$ from the  domain $D^R$ where $u<u_c$. 
A simple coordinate transformation $(x,t)\to(y,t)$ with $y=x-s(t)$ fixes the
boundary at the origin and transforms the
domains $D^L$ and $D^R$ into $Q^L= \mathbb{R}^- \times \mathbb{R}^+$ and $Q^R = \mathbb{R}^+ \times \mathbb{R}^+$
and the moving boundary   problem  becomes 
the following
equivalent initial-boundary value problem that we refer to as QIVP 
(with a detailed derivation given in Part I):
\begin{linenomath}
\begin{subequations} \eelab{QIVP}
\begin{align}
& u_t - \dot{s}(t) u_y  = u_{yy} + f_c(u), \quad (y,t) \in Q^L \cup Q^R ,    \eelab{QIVPa}  \\
& u  \geq u_c \mbox{  in  } \bar{Q}^L, \quad u \leq u_c \mbox{  in  } \bar{Q}^R,      \\
& u(y,0)  = \begin{cases}  1, \quad & y < 0  \\
 0, \quad  & y \geq 0  \end{cases} \eelab{QIVPb} \\
& u(y,t)  \to \begin{cases}  1, \quad &\mbox{as} \quad y \to - \infty  \\
 0, \quad & \mbox{as} \quad y \to \infty \end{cases}   \eelab{QIVPc}
\end{align}
uniformly for $t \in [0,T]$ for all $T >0$. At the boundary,
\begin{align}
& u(0,t)  = u_c, \quad t \in (0,\infty),  \eelab{QIVPd}  \\
& u_y(0^+,t)  = u_y(0^-,t), \quad t \in (0,\infty).  \eelab{QIVPe} \\
& s(0^+)  = 0. \eelab{QIVPf}
\end{align} 
\end{subequations}
\end{linenomath}
In Part I we  stated   regularity conditions  (see equation (18)) 
for the solution $u(y,t)$  and $s(t)$ to be classical for all $t>0$, and on using the classical 
	  maximum principle and comparison theorem (see, for example, \cite{AronsonSerrin1967} and \cite{Fife1979}), obtained that 
	  $0 < u(y,t) < u_c$ for all  $(y,t) \in Q^R$,   
	  $u_c < u(y,t) < 1$ for all $(y,t) \in Q^L$,
	   and  $u_y(y,t)<0$ for all $t > 0$ and $y\in\mathbb{R}$ with $[u_{yy} (y,t)]^{y = 0^+}_{y = 0^-} = f_c^+$ for all  $t \in \mathbb{R}^+$ with $f_c^+ = f_c(u_c^+)$.
We then established that 
in the presence of a cut-off, the initial-boundary value problem \eeref{QIVP}
admits   exactly one  PTW   solution (up to translation) $u(y,t)=U_T(y)$ 
that is   strictly monotone decreasing and positive,  
with $\lim_{y\to-\infty}U_T(y)=1$ 
and $\lim_{y\to\infty}U_T(y)=0$
where  the speed $v=v^*(u_c)$  is, 
for fixed $u_c\in(0,1)$, uniquely defined. 
An explicit expression of $v^*(u_c)$ is in general not known, it is however straightforward to establish that  $v^*(u_c)$ is a
 continuous, monotone decreasing function of $u_c\in(0,1)$, with 
$v^*(u_c)\to 2^-$ as $u_c\to 0^{+}$
and $v^*(u_c)\to 0^+$ as $u_c\to 1^{-}$  \cite{Tisbury_etal}.
Brunet and Derrida \cite{BrunetDerrida1997} predicted
that the difference between  $v^*(u_c)$ and $v_m=2$   
 is  strongly influenced at  small values of $u_c$,
being only logarithmically small in $u_c$ as 
$u_c\to 0^+$. 
This behaviour was rigorously verified by Dumortier, Popovic and Kaper
 \cite{Dumortier_etal2007}, 
 with higher order corrections obtained in Part I. This behaviour is in contrast with the behaviour of  
$v^*(u_c)$ obtained as $u_c\to 1^-$  in which case it vanishes algebraically  in $(1-u_c)$ 
(see Part I).

    \begin{figure}
		 \begin{center}
 	     \includegraphics[width=0.95\textwidth]{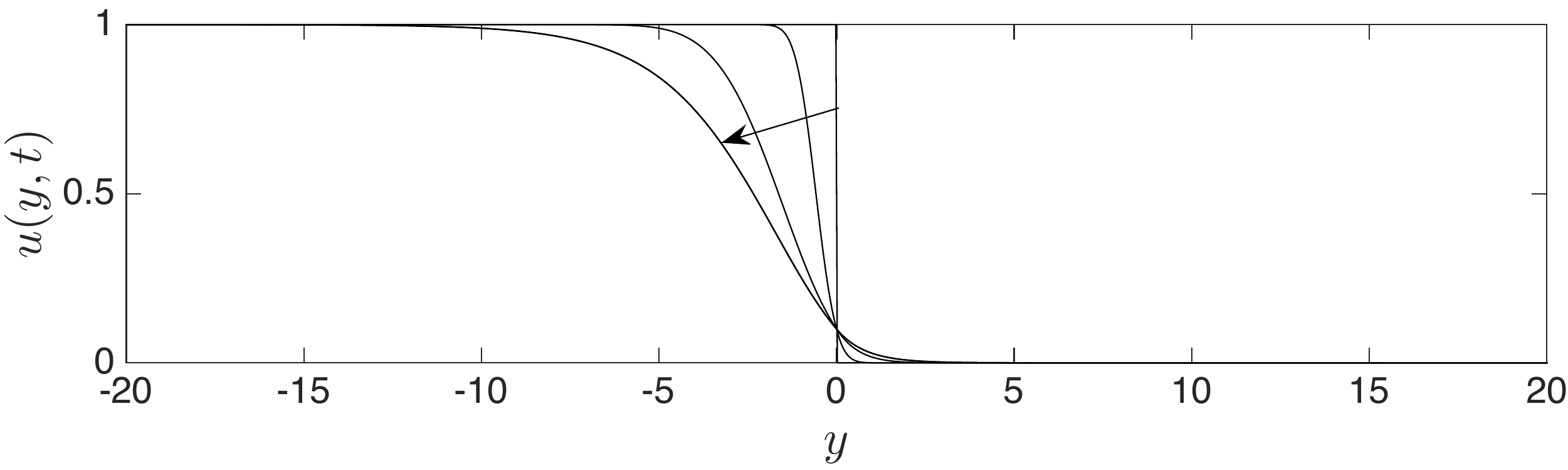}\\
		 (a) $u_c =0.1$
 	     \includegraphics[width=0.95\textwidth]{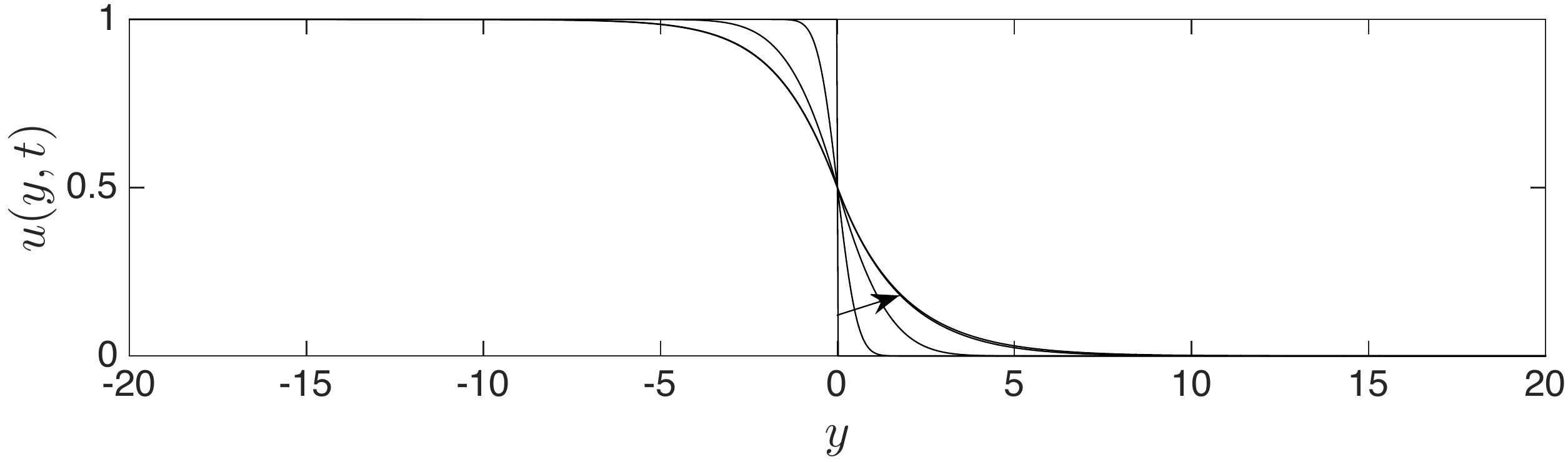}\\
		 (b) $u_c =0.5$
 	     \includegraphics[width=0.95\textwidth]{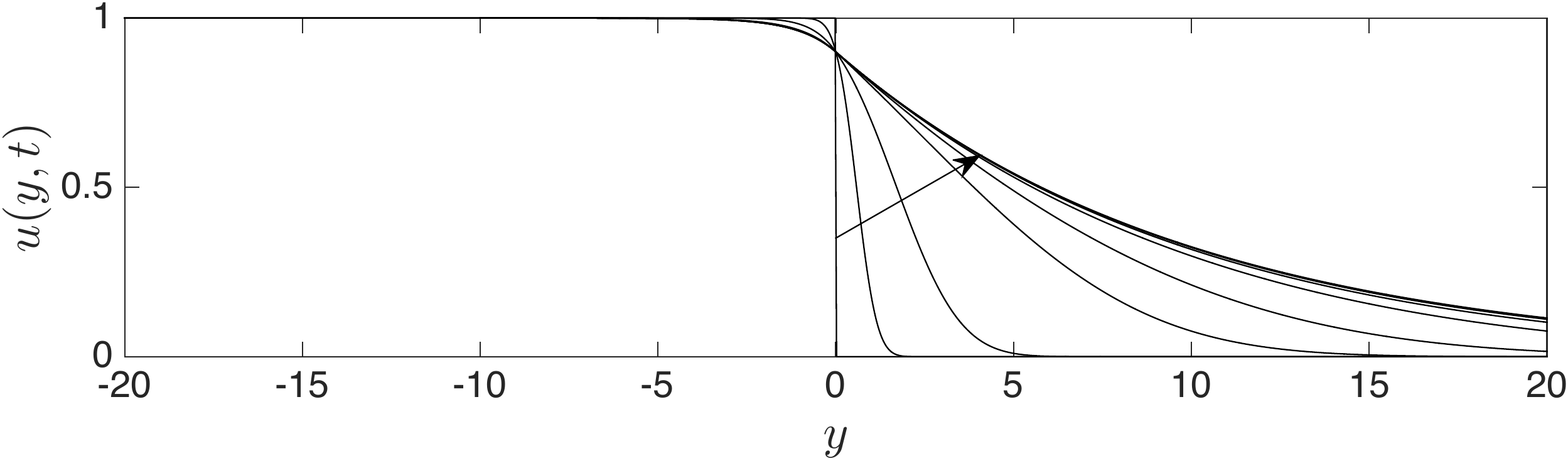}\\
		 (c) $u_c =0.9$
		  \end{center}
 \caption{A graph of the solution $u(y, t)$ to QIVP  as it evolves over time. Results are obtained numerically for (a) $u_c=0.1$, (b) $u_c=0.5$ and (c) $u_c=0.9$ for  $t=0$, $0.1$, $1$, $10$ and $t=30$ with the arrow pointing in the direction of increasing $t$. Panel (c) includes  additional graphs of solutions obtained at $t=100$, $200$, $300$, $350$ and $t=400$. }
  \fflab{PTWu}  
  \end{figure}

We may now once again enquire as to whether or not  a PTW  solution
evolves in the solution to \eeref{QIVP} for arbitrary cut-off $u_c \in (0,1)$  at large time, and, if this is the case, what is the rate of convergence onto the PTW  solution. 
In this paper we  observe
that a PTW  of speed $\lim_{t\to\infty}\dot s (t)=v^*(u_c)$ emerges in 
the solution of \eeref{QIVP} for $t\to\infty$  via numerical simulations obtained for the specific case of $f_c$ with $f$ given by
\eeref{Fisher}.  
We then adapt the approach introduced in \cite{LeachNeedham2003}, where $u_c=0$, to obtain
the detailed description of the large-$t$ structure of the solution to \eeref{QIVP}. 
In particular, we use the theory of matched asymptotic coordinate expansions to 
establish that for  each value of $u_c\in(0,1)$, the solution to  \eeref{QIVP} 
 converges to
the PTW  solution  with propagation speed $v=v^*(u_c)$
at a rate that is linearly exponentially small  in $t$
as $t\to\infty$, specifically
 $O(\dot s(t)-v^*(u_c))$, 
where 
\beq\eelab{expo_convergence}
\dot s(t)=v^*(u_c)+O\left(t^\gamma\exp\left(-\frac{1}{4}v^*(u_c)^2t\right)\right),\quad\text{as $t\to\infty$},
\eeq
(with $\gamma=-1/2$ or $-3/2$ depending on the structure  of $f(u)$, 
specifically   $f'(U_T)$, 
which determines the solution to \eeref{BVP_base_fun_5l}  on which  
the choice in the value of $\gamma$ depends) 
  so that convergence slows down as $u_c$ increases. 
Thus, introducing an arbitrary cut-off   into the reaction function 
changes  the rate of convergence of the large-time solution onto the PTW from algebraic to exponential.  The paper is organised as follows: in section 2, we present numerical results for the specific case of the cut-off Fisher reaction function with $f$ given by
\eeref{Fisher}. Sections 3 and 4 
are respectively devoted to  the small-$t$ ($y\in\mathbb{R}$) and
  large-$|y|$ ($t\geq O(1)$)  structure of the
solution to QIVP. 
These are used in section 5 to develop the complete  asymptotic structure to QIVP as $t\to\infty$,
uniformly in $y\in\mathbb{R}$.
At the end of sections 3 and  5, we illustrate the theory  
for the specific case of the cut-off Fisher reaction function (for which $\gamma=-3/2$).
The paper ends with the concluding section 6.

 \begin{figure} 
\centering
{\includegraphics[width=0.55\linewidth]{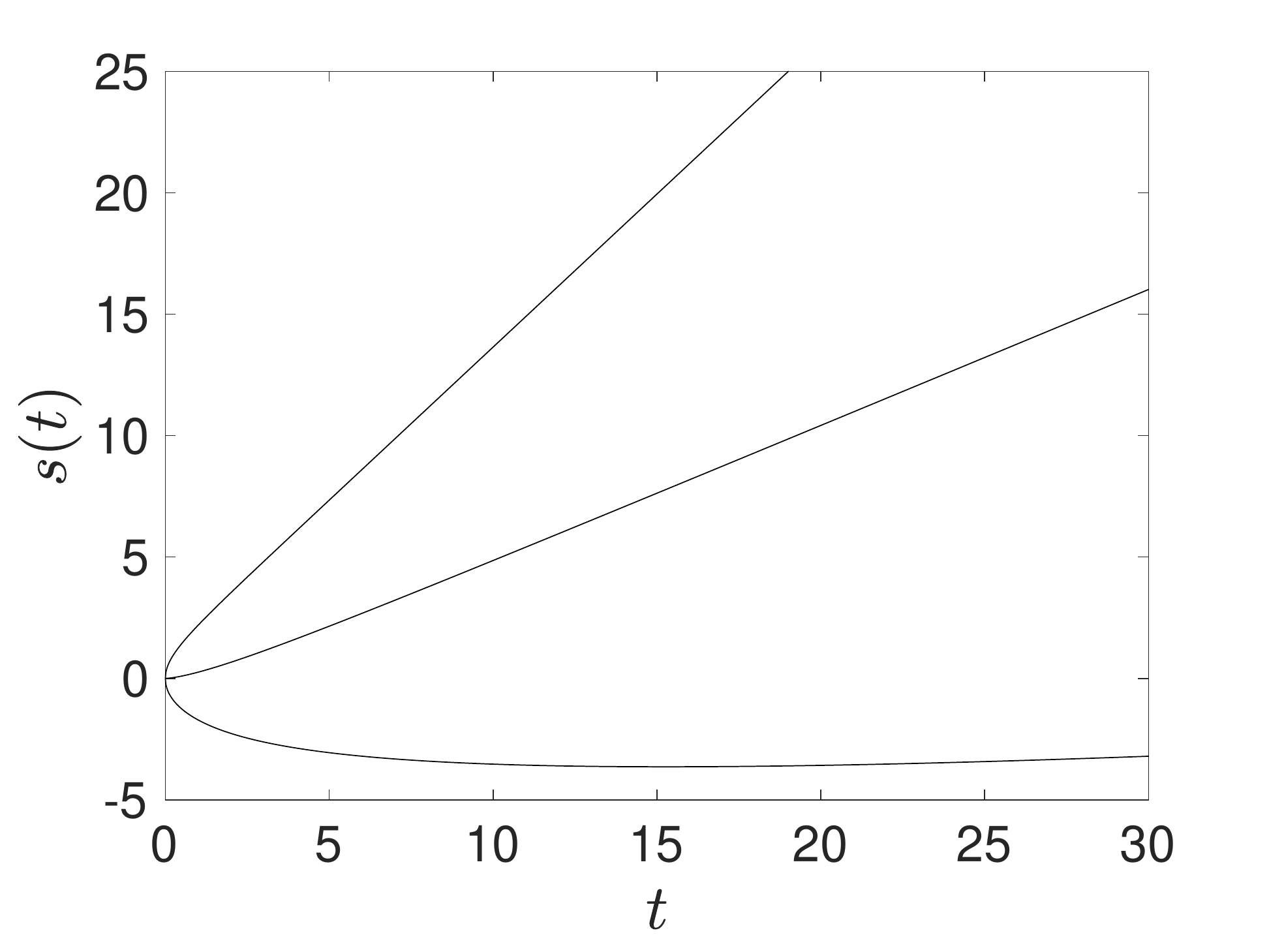}} 
 \caption{A graph of the solution $s(t)$ to QIVP obtained numerically for $u_c=0.1$ (top), $u_c=0.5$ (middle) and $u_c=0.9$ (bottom). }
 \fflab{PTWs}  
 \end{figure}

\section{Numerical solution to QIVP}   \sslab{Numerical_results}

 In this section we consider a numerical solution to QIVP to indicate whether the solution converges onto a PTW solution at large times. 
 We present results for the particular case of the cut-off Fisher reaction function, namely,
 \beq \eelab{cutoff_Fisher}
 f_c(u) = \begin{cases}  u(1-u), \quad & u \in (u_c, \infty) ,     \\
   0, \quad &  u \in (- \infty, u_c],  \end{cases}  
 \eeq
 for fixed cut-off value $u_c \in (0,1)$.
  We adopt an explicit  finite difference scheme, detailed in  Appendix \ssref{Numerical scheme}. 
  We choose this scheme over an implicit scheme despite the severe numerical  stability restrictions on the time step. This is because an explicit scheme is very straightforward to use: at each time step, the  associated numerical calculation
  requires the solution of a linear algebraic system 
  (rather than a nonlinear algebraic system that would be required for an implicit scheme).

We examine the  behaviour of $u(y,t)$, $s(t)$ and $\dot{s}(t)$, obtained numerically for illustrative values of $u_c \in (0,1)$. 
Figures \ffref{PTWu}--\ffref{sdot1}  respectively  focus on the structure of $u(y,t)$, $s(t)$ and  $\dot{s}(t)$ obtained for $u_c=0.1$, $0.5$ and $0.9$.
These confirm all of the qualitative properties obtained in Part I (see equation (20)) and described in section 1. 
 \begin{figure} [t] 
 		     \centering
 		     \begin{minipage}{0.485\textwidth}
 		  	     \includegraphics[width=\textwidth]{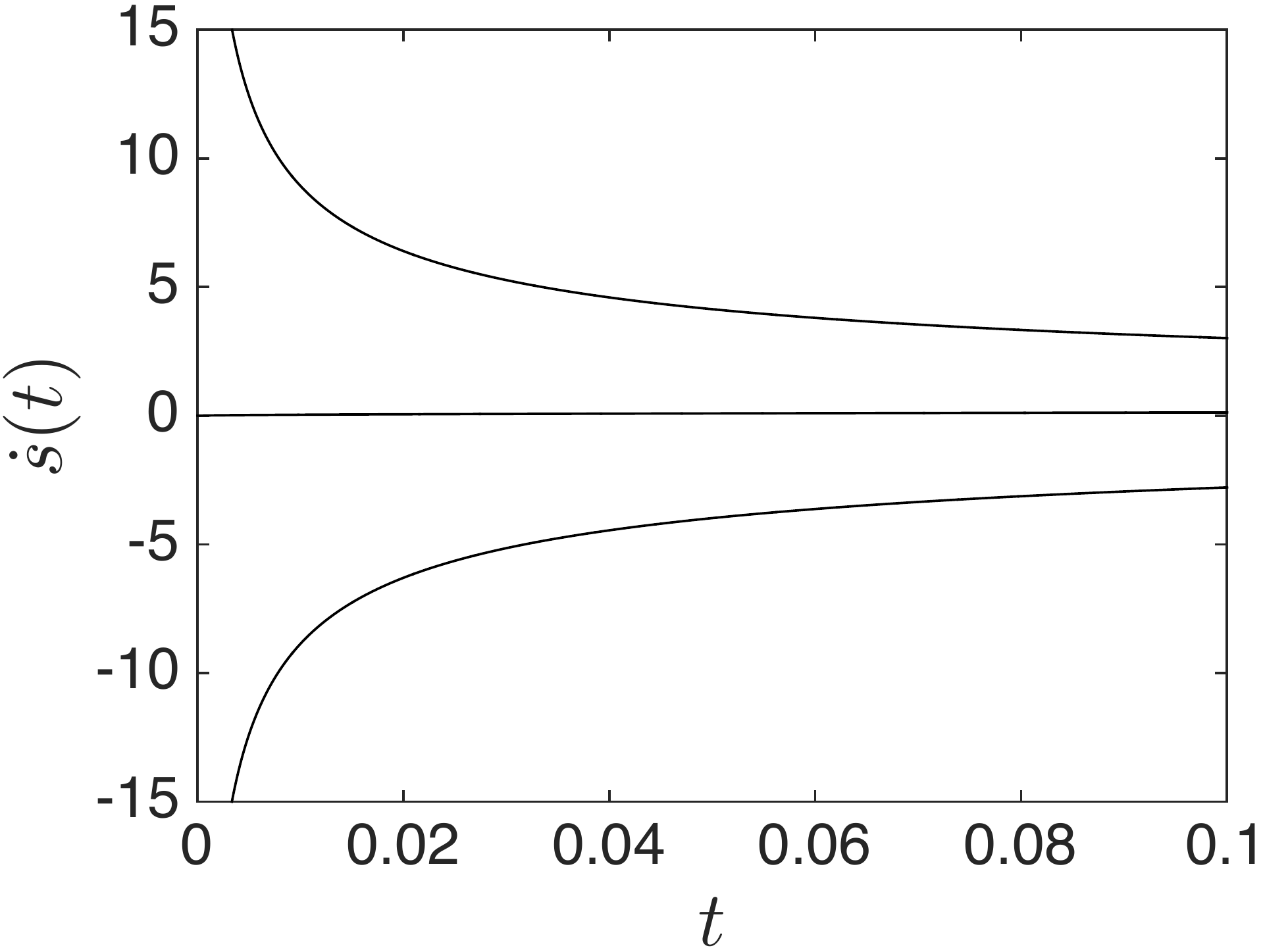}\\
 		  	      \centering  (a) 
 		     \end{minipage}
 		     \hfill
 		     \begin{minipage}{0.485\textwidth}
 		  	   \includegraphics[width=\textwidth]{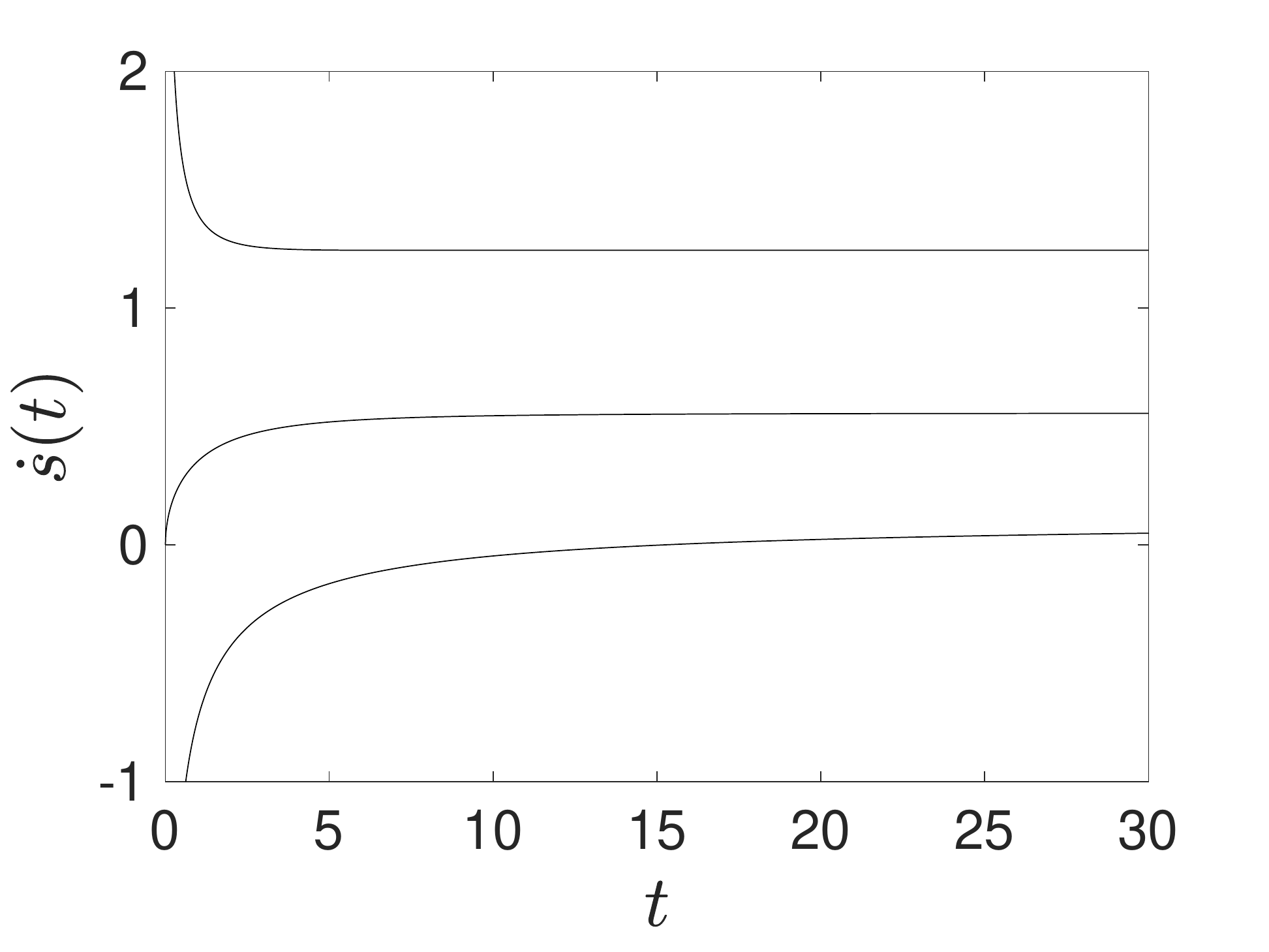}\\
 		    \centering  (b)  
			   \end{minipage}
 \caption{A graph of  $\dot{s}(t)$ to QIVP (solid lines)  obtained numerically for cut-off  value  $u_c=0.1$ (top), $0.5$ (middle) and $0.9$ (bottom)  
 plotted  for a  (a) small and (b) large  range of values of $t$.
 }
 \fflab{sdot1}
 \end{figure} 
  \begin{figure}[h]
  		     \centering
  		     \begin{minipage}{0.485\textwidth}
  		  	     \includegraphics[width=\textwidth]{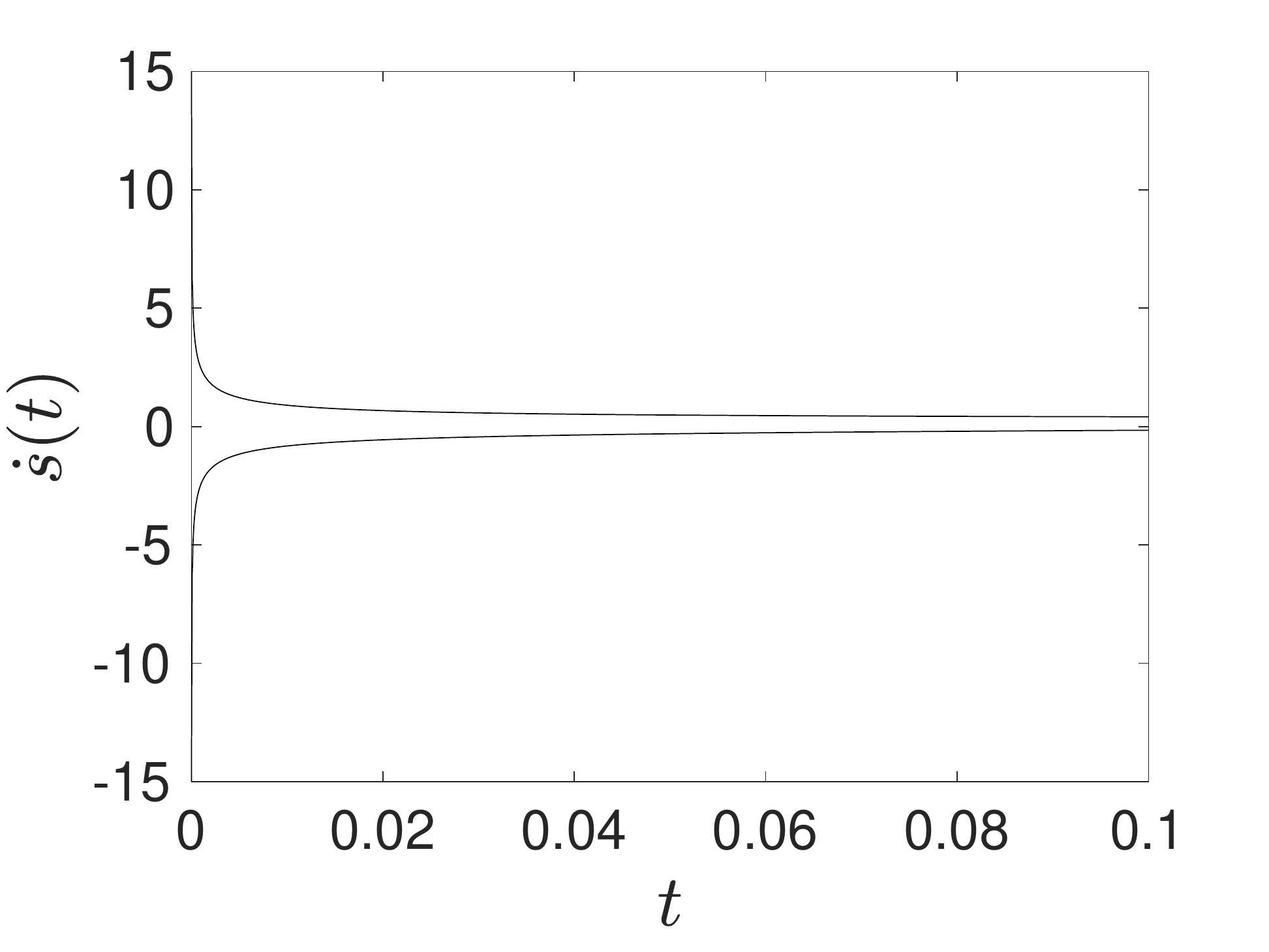}\\
  		  	     \centering  (a)
  		     \end{minipage}
  		     \hfill
  		     \begin{minipage}{0.485\textwidth}
  		  	   \includegraphics[width=\textwidth]{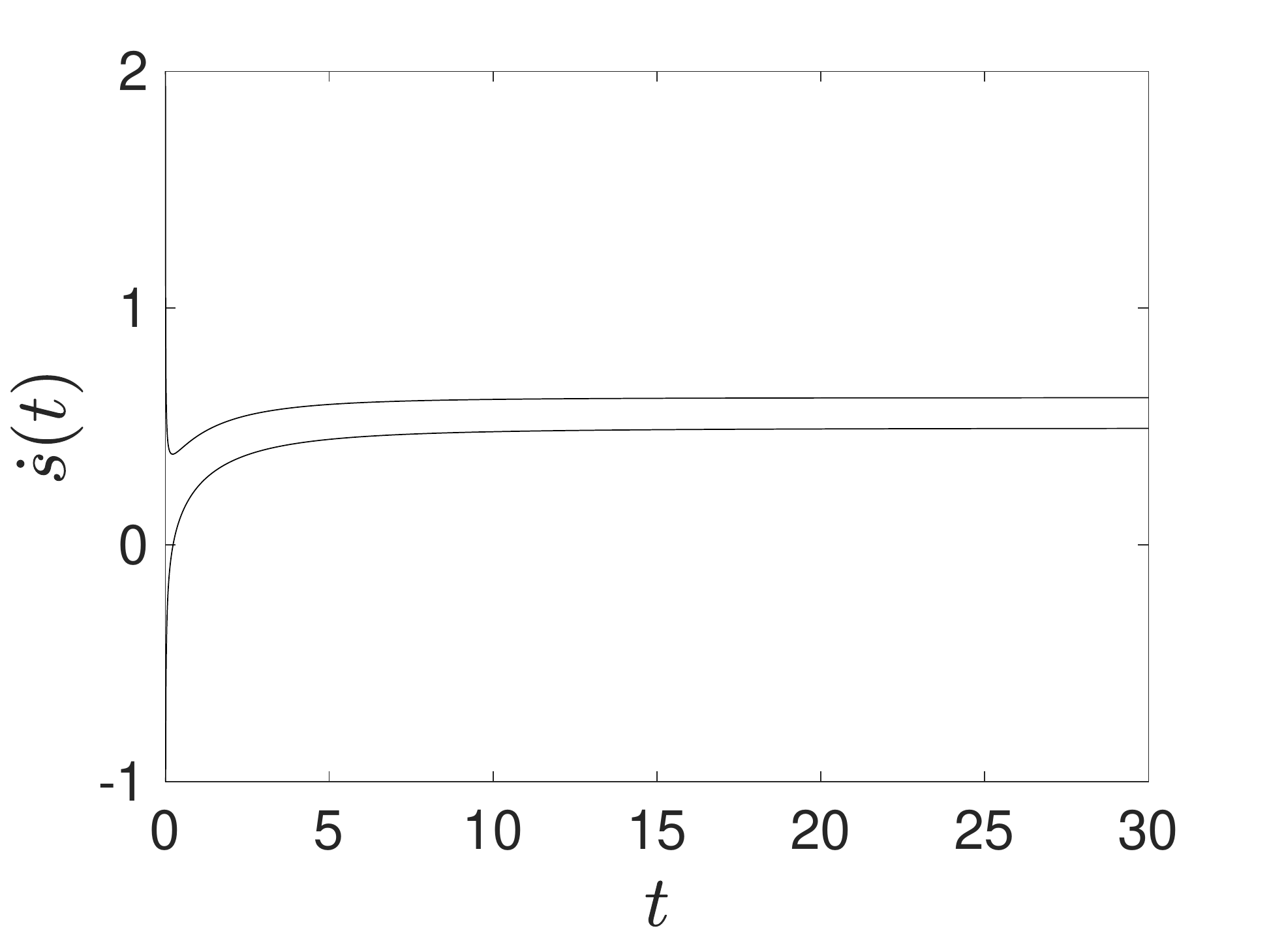}\\
  		    \centering  (b)
  		     \end{minipage}
 \caption{
 Same as Figure \ffref{sdot1} but this time  $u_c=0.45$ (top) and $u_c=0.55$ (bottom).
 }
  \fflab{sdot2}
  \end{figure}
Figure \ffref{PTWu}   indicates that a PTW develops in the large-time structure of the solution to QIVP, that is, as $t\to\infty$. 
   Moreover, the rate of convergence of the solution to the PTW depends on the value of $u_c$ (compare panel (a) with panel (c)).  
Figures \ffref{PTWs} and \ffref{sdot1} show that this PTW will  have propagation speed given by $\lim_{t \to \infty} \dot{s}(t)= v_{\infty}(u_c)$
and in this case, this limit has 
\beq
v_{\infty}(u_c)\simeq \begin{cases}
		 1.248, & \text{for $u_c=0.1$},\\
	             0.558, & \text{for $u_c = 0.5$}, \\
	             0.100, & \text{for $u_c = 0.9$}.
	\end{cases}
\eeq
Figure \ffref{sdot1} also illustrates that $\dot{s}(t)$ appears to have a (integrable) singularity at $t = 0^+$ when $u_c\neq 0.5$. 
This is further supported in Figure \ffref{sdot2} which shows the behaviour of $\dot{s}(t)$  when $u_c=0.45$ and $u_c=0.55$.
For $u_c=0.5$, Figure \ffref{sdot1}  suggests that $\dot{s}(t)$ is regular in this limit, tending to $0$ from above.  
Figures \ffref{sdot1} and  \ffref{sdot2} show that the sign of $\dot{s}(t)$ as $t \to 0^+$ depends upon $u_c$, with $\dot{s}(t)$ initially positive when $0 < u_c < 0.5$ and initially negative when $0.5 < u_c < 1$. Moreover, when $0 < u_c \lesssim 0.2$, then $\dot{s}(t)$ is monotonic decreasing for all $t > 0$; when 
$0.2 \lesssim u_c < 0.5$, then $\dot{s}(t)$ decreases to a minimum value, before increasing to $v_{\infty}(u_c)$; and when $0.5<u_c < 1$, then $\dot{s}(t)$ is monotonic increasing for all $t > 0$. 
 Finally, the correction to $\dot s(t)$ as $t\to\infty$ appears to be exponentially   small in $t$. 
 These features are persistent for all considered values of $u_c\in(0,1)$.

 \begin{figure}[t]
\centering
 {\includegraphics[width=0.55\linewidth]{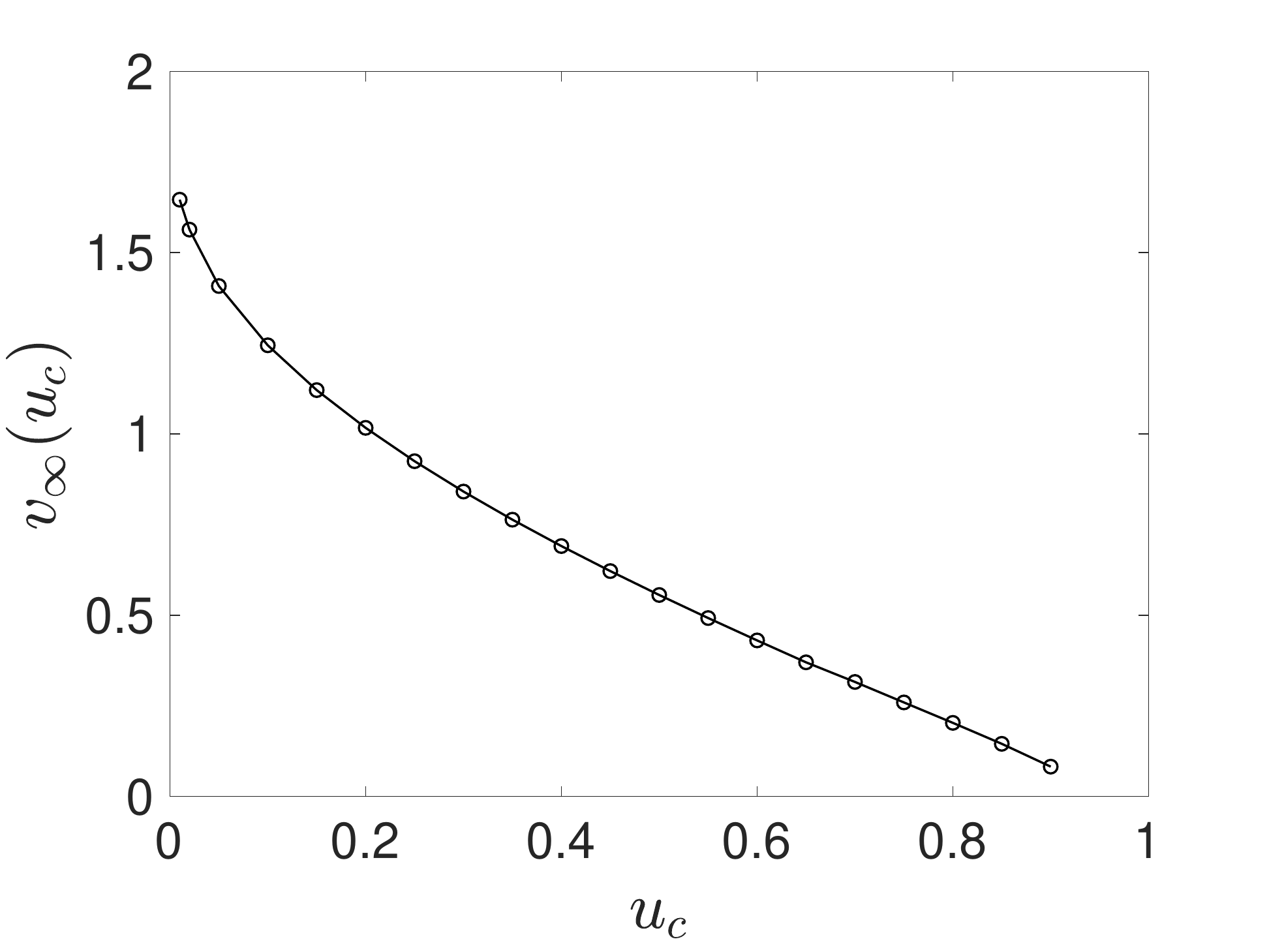}}
 \caption{A graph of   $\lim_{t \to \infty} \dot{s}(t)= v_{\infty}(u_c)$ 
 obtained from the numerical solution to QIVP 
 for selected values of $u_c \in (0,1)$. 
 }
 \fflab{vinf}
 \end{figure}

We conclude that the numerical solution of QIVP involves the formation of a PTW as $t \to \infty$, which has propagation speed  $v_{\infty}(u_c)$  for all values of $u_c \in (0,1)$. 
A graph of numerically calculated values $v_{\infty}(u_c)$ for $u_c \in (0,1)$ is given in Figure \ffref{vinf}, which indicates that $v_{\infty}(u_c)$ is monotone decreasing with $u_c \in (0,1)$. The numerical cost increases drastically as $u_c\to 0^+$ and $u_c\to 1^-$. Nevertheless, we expect that $v_{\infty}(u_c) \to 2^-$ as $u_c \to 0^+$, whilst, $v_{\infty}(u_c) \to 0^+$ as $u_c \to 1^-$. 
Finally, it is instructive to 
compare the travelling wave speed obtained in the large-time limit of the numerical solution to QIVP, namely $v_{\infty}(u_c)$, with a permanent form travelling wave propagation speed, $v^*(u_c)$, obtained numerically in Part I.
As anticipated, we find that, with a significant degree of accuracy (at least up to  two decimal places),  $v_{\infty}(u_c)\approx v^*(u_c)$.


\section{Asymptotic solution  to QIVP as $t \to 0^+$} \sslab{small_time}


We now develop the asymptotic structure to QIVP as $t \to 0^+$ via the method of matched asymptotic coordinate expansions. We anticipate that the structure of the solution to QIVP as $t \to 0^+$ will have two asymptotic regions in $y<0$, and two asymptotic regions in $y>0$. An examination of the leading order balances in equation \eeref{QIVPa}, together with the initial condition \eeref{QIVPb} and the connection conditions \eeref{QIVPd}, \eeref{QIVPe} determine the asymptotic structure as: 
\begin{subequations} \eelab{small_t_structure}
 \begin{align}
& \mbox{region $\mathbf{I_L}$}:  \hspace{0.3cm}  y = O(t^{\frac{1}{2}}) < 0 \mbox{ with }  u=O(1) \mbox{  as  } t \to 0^+,  \eelab{region_1l}  \\
& \mbox{region $\mathbf{I_R}$}:   \hspace{0.23cm} y = O(t^{\frac{1}{2}}) > 0 \mbox{ with }  u=O(1) \mbox{  as  } t \to 0^+, \eelab{region_1r} \\
& \mbox{region $\mathbf{II_L}$}:  \hspace{0.1cm}  y = O(1) < 0 \mbox{ with }  u=1+o(1) \mbox{  as  } t \to 0^+, \eelab{region_2l}  \\
& \mbox{region $\mathbf{II_R}$}:  \hspace{0.07cm}  y = O(1) > 0 \mbox{ with }  u=o(1) \mbox{  as  } t \to 0^+.   \eelab{region_2r}
\end{align}
\end{subequations}
 The situation is illustrated in Figure \ffref{small_structure} (for any variable $\lambda$, we will henceforth write $\lambda=O(1)>0$ as $\lambda=O(1)^+$, and correspondingly, $\lambda=O(1)<0$ as $\lambda=O(1)^-$). It follows from the small-time asymptotic structure \eeref{small_t_structure} of QIVP that we anticipate an asymptotic expansion for $s(t)$ of the form 
\beq \eelab{exp_small_s}
s(t) = s_0 t^{\alpha}  + s_1 t^{\beta} + o(t^{\beta})  \quad \mbox{as} \quad  t \to 0^+, 
\eeq
where the constants $s_0$, $s_1$, $\alpha$ and $\beta(>\alpha)$ are to be found. The initial condition \eeref{QIVPf}, together with a leading order  balance in equation \eeref{QIVPa} determines
\beq
\alpha=\frac{1}{2}.
\eeq
\subsection{Regions  $\mathbf{I_L}$ and $\mathbf{I_R}$ }
We begin in region $\mathbf{I_L}$, following \eeref{region_1l}, where we introduce the coordinate $\eta = y t^{- \frac{1}{2}}=O(1)^-$ as $t \to 0^+$ and where 
$u=u(\eta,t)$ satisfies, from  \eeref{QIVPa},
\beq \eelab{PDE_eta1}
u_t  - \frac{1}{t} \frac{\eta}{2}  u_{\eta} -  \frac{\dot{s}(t)}{t^{\frac{1}{2}}} u_{\eta}  = \frac{1}{t} u_{\eta \eta} + f(u), \quad \eta < 0.
\eeq
We expand $u(\eta,t)$ in the form, 
\beq \eelab{exp_u_1l}
u(\eta, t) = u_{L0} (\eta) + \phi_L(t) u_{L1}(\eta) + o(\phi_L(t)) \quad \mbox{as} \quad  t \to 0^+, 
\eeq
with $\eta =O(1)^-$ and $\phi_L(t)=o(1)$ as $t \to 0^+$ to be determined. 
On substituting  expansions \eeref{exp_small_s} and \eeref{exp_u_1l} into equation \eeref{PDE_eta1}, we obtain at leading order as $t \to 0^+$, 
\begin{subequations} \eelab{BVP1l}
\begin{linenomath}
\beq
u_{L0}'' + \frac{1}{2} (\eta + s_0) u_{L0}' = 0, \quad \eta < 0, \eelab{BVP1la}
\eeq
which must be solved subject to the boundary condition \eeref{QIVPd} at $\eta=0$, together with the matching condition with region $\mathbf{II_L}$ as $\eta \to - \infty$. Using \eeref{region_2l} and \eeref{exp_u_1l}, these conditions require,
\begin{figure}
\centering
{\includegraphics[width=0.65\linewidth]{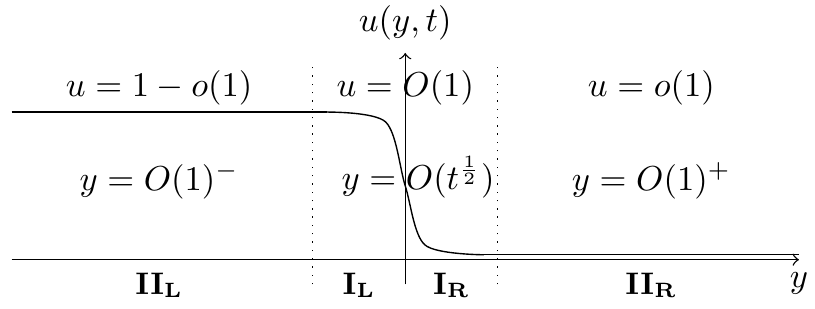}} 
\caption{A sketch of the structure of the solution to QIVP as $t \to 0^+$.}
\fflab{small_structure}
\end{figure}
\begin{align} 
u_{L0} (0) & = u_c,  \\
u_{L0} (\eta) & \to 1 \quad \mbox{as} \quad \eta \to - \infty .
\end{align}
\end{linenomath}
\end{subequations}
Due to the coupling condition \eeref{QIVPe} across $y=0$, it is necessary now to consider region $\mathbf{I_R}$, in which, via \eeref{region_1r}, $\eta =O(1)^+$ and $u=O(1)$ as $t \to 0^+$ and where 
$u=u(\eta,t)$ satisfies, from  \eeref{QIVPa}, 
\beq \eelab{PDE_eta2}
u_t  - \frac{1}{t} \frac{\eta}{2}  u_{\eta} -  \frac{\dot{s}(t)}{t^{\frac{1}{2}}} u_{\eta}  = \frac{1}{t} u_{\eta \eta}, \quad \eta > 0.
\eeq
We expand $u(\eta,t)$ in the form, 
\beq \eelab{exp_u_1r}
u(\eta, t) = u_{R0} (\eta) + \phi_R(t) u_{R1}(\eta) + o(\phi_R(t)) \quad \mbox{as} \quad  t \to 0^+,
\eeq
with $\eta =O(1)^+$ as $t \to 0^+$. Here $\phi_R=o(1)$ as $t \to 0^+$, and is to be determined. Now, substituting  expansions \eeref{exp_small_s} and \eeref{exp_u_1r} into equation \eeref{PDE_eta2}, we obtain at leading order as $t \to 0^+$,  
\begin{subequations} \eelab{BVP1r}
\beq
u_{R0}'' + \frac{1}{2} (\eta + s_0) u_{R0}' = 0, \quad \eta > 0, \eelab{eq:LO equation in I R}
\eeq
which must be solved subject to the boundary condition \eeref{QIVPd} at $\eta=0$, together with the matching condition with region $\mathbf{II_R}$ as $\eta \to  \infty$, which requires,
 \begin{align} 
u_{R0} (0) & = u_c,  \\
u_{R0} (\eta) & \to 0 \quad \mbox{as} \quad \eta \to  \infty .
\end{align}
\end{subequations}
Finally, the boundary value problems \eeref{BVP1l} and \eeref{BVP1r} must be solved subject to the coupling condition \eeref{QIVPe} across $\eta=0$, which requires
\beq \eelab{BVP1c}
u_{L0}' (0) = u_{R0}' (0).
\eeq
The solutions to \eeref{BVP1l} and \eeref{BVP1r} respectively, are readily obtained as
\begin{subequations} \eelab{gen_sol1}
 \begin{align} 
u_{L0}(\eta) & = \frac{ u_c \left( 1 +  \erf(\frac{\eta + s_0}{2}) \right)  -  \erf(\frac{\eta + s_0}{2}) + \erf(\frac{s_0}{2}) }{ \left( 1+ \erf(\frac{s_0}{2}) \right) }  , \quad \eta \leq 0,     \eelab{gen_sol1l}   \\
u_{R0}(\eta) & =  u_c  \frac{ \left( 1 - \erf(\frac{\eta + s_0}{2}) \right) }{ \left( 1- \erf(\frac{s_0}{2}) \right) } , \quad \eta \geq 0. \eelab{gen_sol1r}
\end{align}
\end{subequations}
Finally, an application of condition \eeref{BVP1c} to \eeref{gen_sol1} determines
\beq \eelab{small_s0}
s_0 = 2 \erf^{-1} (1- 2 u_c),
\eeq
and thus, the leading order terms in region $\mathbf{I_L}$ and region $\mathbf{I_R}$, respectively, are given by
\begin{subequations} \eelab{par_sol1}
 \begin{align} 
u_{L0}(\eta) & =  \frac{1}{2} \bigg[ 1 - \erf \left( \frac{\eta}{2} + \erf^{-1} (1 - 2 u_c) \right) \bigg],  \quad \eta \leq 0, \eelab{par_sol1l} \\
u_{R0}(\eta) & =   \frac{1}{2} \bigg[ 1 - \erf \left( \frac{\eta}{2} + \erf^{-1} (1 - 2 u_c) \right) \bigg] , \quad \eta \geq 0 \eelab{par_sol1r} . 
\end{align}
\end{subequations}

We now proceed to the correction terms in expansions \eeref{exp_small_s}, \eeref{exp_u_1l} and \eeref{exp_u_1r}. A balancing of terms requires $\phi_L(t)= \phi_R(t)=O(t)$ as $t \to 0^+$ and $\beta = \frac{3}{2}$. Thus, we set $\phi_L(t)= \phi_R(t)=t$, without loss of generality. On substitution from expansions \eeref{exp_small_s}, \eeref{exp_u_1l} and \eeref{exp_u_1r} into equations \eeref{PDE_eta1} and \eeref{PDE_eta2}, we obtain the coupled problem for $u_{L1}(\eta) (\eta<0)$, $u_{R1}(\eta) (\eta>0)$ and $s_1$, namely, 
\begin{subequations} \eelab{2nd_BVP1}
\begin{align}
& u_{L1}'' + \frac{1}{2} (\eta + s_0) u_{L1}' - u_{L1} = -\frac{3}{2} s_1 u_{L0}' - f(u_{L0}(\eta)) ,  \quad \eta < 0,  \eelab{2nd_BVP1l}  \\
& u_{R1}'' +  \frac{1}{2} (\eta + s_0) u_{R1}' - u_{R1}  =  -\frac{3}{2} s_1 u_{R0}', \quad \eta > 0,  \eelab{2nd_BVP1r}
\end{align} 
subject to the coupling conditions 
\begin{align} 
u_{L1} (0) & = u_{R1}(0) = 0,  \eelab{2nd_BVPc1}  \\
u_{L1}' (0) & = u_{R1}'(0),  \eelab{2nd_BVPc2}
\end{align} 
and the matching conditions to region $\mathbf{II_L}$ and to region $\mathbf{II_R}$, respectively, which are readily obtained as, 
\begin{align} 
u_{L1} (\eta) & \to 0 \quad \mbox{as} \quad \eta \to - \infty , \eelab{2nd_BVPlb} \\
u_{R1} (\eta) & \to 0 \quad \mbox{as} \quad \eta \to  \infty . \eelab{2nd_BVPrb}
\end{align}
\end{subequations} 
In considering the coupled problem \eeref{2nd_BVP1}, we first observe that $1 + \frac{1}{2} (\eta + s_0)^2$ is a solution to the homogeneous Part of both \eeref{2nd_BVP1l} and \eeref{2nd_BVP1r}. With this observation, together with the method of 
variation of parameters, we can write the general solutions to \eeref{2nd_BVP1l} and \eeref{2nd_BVP1r} as, 
\begin{linenomath}
\begin{subequations}
\begin{align} 
u_{L1}(\eta)  = \; & d_1 \hat{u}(\eta) + d_2 \bar{u}(\eta)  - \frac{s_1}{2 \sqrt{\pi}}  \exp \left( - \left( \frac{\eta + s_0}{2} \right)^2 \right)+ u_{p2}(\eta)  , \quad \eta \leq 0, \eelab{2nd_gen_sol1l} \\
u_{R1}(\eta) =  \; & \bar{d}_1 \hat{u}(\eta)  + \bar{d}_2  \bar{u}(\eta)  - \frac{s_1}{2 \sqrt{\pi}}  \exp \left( - \left( \frac{\eta + s_0}{2} \right)^2 \right) , \quad \eta \geq 0,  \eelab{2nd_gen_sol1r}
\end{align}
\end{subequations}
\end{linenomath}
where $d_1, d_2, \bar{d}_1$ and $\bar{d}_2$ are arbitrary constants to be determined and the function $u_{p2}(\eta)$ is given by
\beq
u_{p2}(\eta) =  \frac{\hat{u}(\eta)}{2}  \int_{\eta}^0 I_1(\lambda) d \lambda - \frac{\bar{u}(\eta)}{2}  \int_{\eta}^0 I_2(\lambda) d \lambda   , \quad \eta \leq 0 ,  \eelab{par_sol}
\eeq
with functions
\begin{linenomath}
\begin{subequations} \eelab{small_t_functions}
\begin{align}
& \hat{u}(\eta)= \sqrt{\pi} \left( 1 + \frac{(\eta + s_0)^2}{2} \right)  \erf \left( \frac{\eta + s_0}{2} \right) +  (\eta + s_0)  \exp \left( - \left( \frac{\eta + s_0}{2} \right)^2 \right),  \\ 
& \bar{u}(\eta) =  1 + \frac{(\eta + s_0)^2}{2}, \\	
& I_1(\eta)= \exp \bigg( \left( \frac{\eta + s_0}{2} \right)^2 \bigg)  \bar{u}(\eta)  f(u_{L0}(\eta)), \\
& I_2(\eta) = \exp \bigg( \left( \frac{\eta + s_0}{2} \right)^2 \bigg)  \hat{u}(\eta)  f(u_{L0}(\eta)).
\end{align}
\end{subequations}
\end{linenomath}
Next, an application of condition \eeref{2nd_BVPc1} requires 
 \begin{align}
& d_2 =  \left( \frac{s_1}{\sqrt{\pi}} -2 d_1 s_0 \right)  \frac{e^{- \frac{{s_0}^2}{4}} }{\left( {s_0}^2+2 \right) }  - d_1 \sqrt{\pi} \erf \left( \frac{s_0}{2} \right) , \eelab{d2}  \\
& \bar{d}_2 =  \left( \frac{s_1}{\sqrt{\pi}} -2 \bar{d}_1 s_0 \right)  \frac{e^{- \frac{{s_0}^2}{4}} }{\left( {s_0}^2+2 \right) }  - \bar{d}_1 \sqrt{\pi} \erf \left( \frac{s_0}{2} \right), \eelab{bar_d2} \
\end{align}
whilst, applying the matching conditions \eeref{2nd_BVPlb} and \eeref{2nd_BVPrb} requires 
\begin{align}
& d_2 = \sqrt{\pi}\left(d_1+\frac{1}{2} \hat{d}_1\right),  \eelab{d1}  \\
& \bar{d}_2 = -  \sqrt{\pi} \bar{d}_1 , \eelab{bar_d1}
\end{align}
with the constant $\hat{d}_1$ given by
\beq \eelab{s1_con}
\hat{d}_1 = \int_{- \infty}^0  \left(\sqrt{\pi} I_1(\lambda)  + I_2(\lambda)\right) d \lambda.
 \eeq
As $u_{p2}'(0)=0$, an application of the coupling condition \eeref{2nd_BVPc2} determines $d_1=\bar{d}_1$ (and thus $d_2=\bar{d}_2$) which finally requires that
\beq 
s_1 = \frac{1}{4}  \left( \sqrt{\pi} ({s_0}^2 +2) \left( 1 - \erf \left( \frac{s_0}{2} \right) \right) e^{\frac{{s_0}^2}{4}}  - 2 s_0  \right)   \hat{d}_1, \eelab{exp_small_s1}
\eeq
 after which  (using\eeref{small_s0}), $d_1$, $\bar{d}_1$, $d_2$, $\bar{d}_2$ follow from   \eeref{d2}, \eeref{bar_d2}, \eeref{d1}   and \eeref{bar_d1}. 

Thus, we have determined that the two-term expansions for $u(\eta, t)$ in region $\mathbf{I_L}$ and region $\mathbf{I_R}$ are given by
\begin{align}
u(\eta,t)  =  & \; \frac{1}{2} \bigg[ 1 - \erf \left( \frac{\eta+ s_0}{2} \right) \bigg] 
\nonumber \\
& +  t \left( d_1 \hat{u}(\eta)   + d_2 \bar{u}(\eta)  -  \frac{s_1}{2 \sqrt{\pi}} \exp \bigg[ - \left( \frac{\eta + s_0}{2} \right)^2 \bigg] + u_{p2}(\eta) \right) + o(t),  \eelab{par_sol_1l_full} \\
\intertext{as $t \to 0^+$ with $\eta = O(1)^-$, and }
u(\eta,t)  =  & \; \frac{1}{2} \bigg[ 1 - \erf \left( \frac{\eta+ s_0}{2} \right)\bigg] + t \left( d_1 \hat{u}(\eta)   + d_2 \bar{u}(\eta)  -  \frac{s_1}{2 \sqrt{\pi}} \exp \bigg[ - \left( \frac{\eta + s_0}{2} \right)^2 \bigg] \right) + o(t),   \eelab{par_sol_1r_full}
\end{align} 
as $t \to 0^+$, with $\eta = O(1)^+$, whilst the
  two-term expansion for $s(t)$ is given by 
\beq
s(t) = s_0 t^{\frac{1}{2}} + s_1 t^{\frac{3}{2}} + o(t^{\frac{3}{2}}) , \eelab{exp_small_s_full}
\eeq
as $t \to 0^+$. Here the constants $d_1$, $d_2$, $s_0$ and $s_1$ are given by \eeref{d1}, \eeref{d2}, \eeref{small_s0} and \eeref{exp_small_s1}, respectively, and the functions $\hat{u}(\eta)$, $\bar{u}(\eta)$, $I_1(\lambda)$, $I_2(\lambda)$ and $u_{p2}(\eta)$ are given by \eeref{small_t_functions} and \eeref{par_sol}, respectively. 
It is worth noting that we have obtained the two term small-time expansions for $s(t)$ without needing to know the precise asymptotic structure of the solution in regions $\mathbf{II_L}$ and $\mathbf{II_R}$. The matching conditions with regions $\mathbf{I_L}$ and $\mathbf{I_R}$, respectively, were sufficient. The asymptotic expansion in regions $\mathbf{II_L}$ and $\mathbf{II_R}$ are now obtained to complete the small-time asymptotic structure.

\subsection{Region  $\mathbf{II_L}$} 

First, from \eeref{par_sol_1l_full} and \eeref{par_sol_1r_full}, we observe that for $(- \eta) \gg 1$, 
\beq \eelab{exp_1l_match}
u(\eta,t) \sim 1 - \frac{1}{\sqrt{\pi}} \frac{1}{ \lvert \eta + s_0 \rvert } \exp \bigg( - \left( \frac{\eta + s_0}{2} \right)^2 \bigg)  ( 1 - O((\eta +s_0)^{-2}) ) ,
\eeq
as $t \to 0^+$, and for $\eta \gg 1$, 
\beq \eelab{exp_1r_match}
u(\eta,t) \sim \frac{1}{\sqrt{\pi}} \frac{1}{ \left( \eta + s_0  \right) } \exp \bigg( - \left( \frac{\eta + s_0}{2} \right)^2 \bigg) ( 1 - O((\eta +s_0)^{-2}) ),
\eeq
as $t \to 0^+$.
Now, as $ \eta \to - \infty$ we move out of region $\mathbf{I_L}$ and into region $\mathbf{II_L}$, in which, via \eeref{region_2l},  $y=O(1)^-$ and $u(y,t)=1 + o(1)$ as $t \to 0^+$. The structure of the expansion in region $\mathbf{I_L}$, for $(- \eta) \gg 1$, (given by \eeref{exp_1l_match}) suggests that in region $\mathbf{II_L}$ we write 
\beq \eelab{exp_u2l}
u(y,t) = 1 - e^{-\frac{H(y,t)}{t}},
\eeq
and expand in the form, 
\beq \eelab{exp_u2la}
H(y,t) = H_0(y) + t^{\frac{1}{2}} H_1(y) + t \ln t H_2(y) + t H_3(y) + o(t),
\eeq
as $t \to 0^+$ with $y=O(1)^-$ and $H_0(y) > 0$ (the $t\ln t$ term arises from the algebraic prefactor of the exponential term in \eeref{exp_1l_match}). We substitute expansions \eeref{exp_u2l} and \eeref{exp_u2la} into equation \eeref{QIVPa} to obtain (on solving at each order of $t$ in turn)
\beq\eelab{gen_sol_2l}
\begin{split}
	 u(y,t) = 1  - \exp \Bigg( & -  \frac{y^2 }{4t}   - \frac{1}{t^{\frac{1}{2}}} \bigg( \frac{s_0}{2} y + D_1 (-y)^{\frac{1}{2}} \bigg)   - D_2 \ln t   \\
	 &  - \bigg( \frac{ \left( 1-2D_2 \right) }{2} \ln  (-y)  + \frac{s_0 D_1}{2} 
	 \frac{1}{ \left(-y \right)^{\frac{1}{2}}} + \frac{{D_1}^2}{4} \frac{1}{ y } + D_3 \bigg) +o(1) \Bigg),
\end{split}
\eeq
as $t \to 0^+$, with $y=O(1)^-$, and where $D_1$, $D_2$ and $D_3$ are arbitrary constants to be determined. It remains to match expansion \eeref{gen_sol_2l} in region $\mathbf{II_L}$ (as $y \to 0^-$) with expansion \eeref{exp_1l_match} in region $\mathbf{I_L}$ (as $\eta \to -\infty$). On applying Van Dyke's matching principle \cite{VanDyke1975}, we readily obtain that 
\beq
  D_1=0, \qquad D_2 = - \frac{1}{2}, \qquad D_3 =\frac{1}{2} \ln \pi + \frac{{s_0}^2}{4}.
\eeq 
Thus, the expansion in region $\mathbf{II_L}$ is given by 
 \beq \eelab{u_2l}
u(y,t) =  1  - \exp \Bigg( -\frac{y^2}{4t}  - \frac{y s_0}{2 t^{\frac{1}{2}}}   + \frac{1}{2} \ln t - \left( \ln ( - y ) + \frac{1}{2} \ln \pi  + \frac{{s_0}^2}{4} \right) +o(1) \Bigg),  
\eeq
as $t \to 0^+$,  with $y=O(1)^-$. Furthermore, we conclude from \eeref{u_2l} that this expansion remains uniform for $(-y) \gg 1$ as $t \to 0^+$.

\subsection{Region $\mathbf{II_R}$} 

Next, as $\eta \to \infty$, we move out of region $\mathbf{I_R}$ and into region $\mathbf{II_R}$, in which, via \eeref{region_2r},  $y=O(1)^+$ and $u(y,t)= o(1)$ as $t \to 0^+$. The structure of the expansion in region $\mathbf{I_R}$, for $\eta \gg 1$, (given by \eeref{exp_1r_match}) suggests that in region $\mathbf{II_R}$ we write 
\beq \eelab{exp_u2r}
u(y,t) = e^{-\frac{\bar{H}(y,t)}{t}},
\eeq
and expand in the form,
\beq \eelab{exp_u2ra}
\bar{H}(y,t) = \bar{H}_0(y) + t^{\frac{1}{2}} \bar{H}_1(y) + t \ln t \bar{H}_2(y) + t \bar{H}_3(y) + o(t),
\eeq
as $t \to 0^+$ with $y=O(1)^+$ and $\bar{H}_0(y) > 0$ 
 (the $t\ln t$ term arises from the algebraic prefactor of the exponential term in \eeref{exp_1r_match}).
Substitution of \eeref{exp_u2r} and \eeref{exp_u2ra} into equation \eeref{QIVPa} gives (on solving at each order of $t$ in turn)
\beq \eelab{gen_sol_2r}
\begin{split}
u(y,t) =  \exp \Bigg( & -   \left( \frac{y^2 }{4t} \right)  - \frac{1}{t^{\frac{1}{2}}} \bigg( \frac{s_0}{2} y + \bar{D}_1   y^{\frac{1}{2}} \bigg)   - \bar{D}_2 \ln t  \\
 &    - \bigg( \frac{ \left( 1-2\bar{D}_2 \right) }{2} \ln  y + \frac{s_0 \bar{D}_1}{2} \frac{1}{ y^{\frac{1}{2}}} + \frac{{\bar{D}_1}^2}{4} \frac{1}{y} + \bar{D}_3 \bigg) +o(1) 
\Bigg),
\end{split}
\eeq
as $t \to 0^+$, with $y=O(1)^+$, and where $\bar{D}_1$, $\bar{D}_2$ and $\bar{D}_3$ are arbitrary constants to be determined. It remains to match expansion \eeref{gen_sol_2r} in region $\mathbf{II_R}$ (as $y \to 0^+$) with expansion \eeref{exp_1r_match} in region $\mathbf{I_R}$ (as $\eta \to \infty$). On applying Van Dyke's matching principle \cite{VanDyke1975}, we readily obtain that 
\beq
 \bar{D}_1=0, \qquad \bar{D}_2 = - \frac{1}{2}, \qquad \bar{D}_3 =\frac{1}{2} \ln \pi + \frac{{s_0}^2}{4}.
\eeq 
Thus, the expansion in region $\mathbf{II_R}$ is given by 
\beq \eelab{u_2r}
u(y,t) = \exp \Bigg( -\frac{y^2}{4t}  - \frac{y s_0}{2 t^{\frac{1}{2}}}   + \frac{1}{2} \ln t  -  \left( \ln y + \frac{1}{2} \ln \pi  + \frac{{s_0}^2}{4}  \right) +o(1) \Bigg),
\eeq
as $t \to 0^+$ and  $y=O(1)^+$. Furthermore, we conclude from \eeref{u_2l} that this expansion remains uniform for $y \gg 1$ as $t \to 0^+$.

\bigskip

The asymptotic structure of the solution to QIVP as $t \to 0^+$ is now complete with the expansions 
 \eeref{u_2l},
\eeref{par_sol_1l_full},
\eeref{par_sol_1r_full}
and
 \eeref{u_2r} 
in regions $\mathbf{II_L}$, $\mathbf{I_L}$, $\mathbf{I_R}$ and $\mathbf{II_R}$.  
 We next use this information to enable us to develop the asymptotic structure of the solution to QIVP as $\lvert y \rvert \to \infty$ with $t=O(1)$. However, before proceeding to this, it is of interest to examine the form of $\dot{s}(t)$ in the small-time limit for all $u_c \in (0,1)$.
It follows from  expression \eeref{exp_small_s_full} that
\beq  \eelab{sdot_small_t}
\dot{s}(t) \sim \frac{1}{2}s_0 t^{- \frac{1}{2}} + \frac{3}{2}s_1 t^{\frac{1}{2}} \quad \mbox{as} \quad t \to 0^+,
\eeq
with $s_0$ and $s_1$ given by equations \eeref{small_s0} and \eeref{exp_small_s1} respectively. 
In particular, we observe from \eeref{small_s0} that $s_0$ is monotonic decreasing in $u_c$ with 
\beq \eelab{s0_limits}
s_0 \to \infty \; \; \mbox{as} \; \; u_c \to 0^+, \quad s_0=0 \; \; \mbox{when} \; \; u_c=\frac{1}{2} \quad \mbox{and} \quad s_0 \to - \infty \; \; \mbox{as} \; \; u_c \to 1^-.
\eeq
Thus, the leading term in \eeref{sdot_small_t} reveals that $\dot{s}(t)$ has an integrable singularity as $t \to 0^+$, with 
\beq \eelab{sdot_small1}
\dot{s}(t) \to + \infty \quad \mbox{as} \quad t \to 0^+,
\eeq
when $0 < u_c < 1/2$, whilst, 
\beq \eelab{sdot_small2}
\dot{s}(t) \to - \infty \quad \mbox{as} \quad t \to 0^+,
\eeq
when $1/2 < u_c < 1$. When $u_c=1/2$,   a transition  occurs with $\dot s(t)$ not   
singular and 
\beq \eelab{sdot_small3}
\dot{s}(t) \to 0  \quad \mbox{as} \quad t \to 0^+. 
\eeq

 \subsection{The  case of a cut-off Fisher reaction}
  We observe that \eeref{sdot_small1}, \eeref{sdot_small2} and \eeref{sdot_small3} agree with the numerical solutions for QIVP obtained for the cut-off Fisher reaction function in section \ssref{Numerical_results}, as illustrated in Figures  \ffref{sdot1} and \ffref{sdot2}.
Moreover,  it is straightforward to establish (via
\eeref{s1_con} and \eeref{exp_small_s1})
that for $u_c=1/2$,
$s_1=s_1^*>0$.
Therefore $\dot{s}(t) \to 0^+$ as $t\to 0^+$.
In addition, it is interesting to note from expression \eeref{sdot_small_t} that when
$u_c$ is close to $1/2$   a local minimum point in the graph of $\dot s(t)$ against $t$
bifurcates singularly from $t=0$ as $u_c$ decreases through $u_c=1/2$.
In particular, the local minimum point when $u_c<1/2$   
is located when
$t=t_m\sim \frac{1}{3}s_0/s_1>0$. 
As $u_c\to \frac{1}{2}^-$, 
  $ \frac{1}{3}s_0/s_1\sim \frac{2}{3}\sqrt{\pi}(1-2u_c)/s_1^*+O((1-2u_c)^2)$ 
 where $s_1^*\simeq 0.28$ is approximated numerically using  \eeref{s1_con} and \eeref{exp_small_s1}. 
 The location of the minimum point increases as 
  $u_c $ decreases,   until $u_c\approx 0.2$ when
 $t_m$ is no longer small and in fact  the local minimum point ceases to exist at this sufficiently low value of $u_c$.
This is also in agreement with the numerical solution of section 2 and in particular Figures \ffref{sdot1} and    \ffref{sdot2}.
A comparison of $\dot s(t)$ and $u(y,t)$ as computed from
 \eeref{par_sol_1l_full}, \eeref{par_sol_1r_full}, \eeref{u_2l} and \eeref{u_2r} with the full numerical solution to QIVP obtained for the cut-off Fisher reaction function
is readily made
(but for brevity is not presented here). This demonstrates
the full agreement
with the small-time asymptotic structure of the solution
obtained in this section and the numerical solution obtained in section \ssref{Numerical_results} for $t$ small.

\section{Asymptotic solution to QIVP  as $\lvert y \rvert \to \infty$  with $t=O(1)$} \sslab{large_spatial}

We now develop the structure of the solution to QIVP as $\lvert y \rvert \to \infty$ with $t=O(1)$. 
\subsection{Region  $\mathbf{III_L}$} \sslab{large_spatiala}
We begin in region $\mathbf{III_L}$, where $y \to - \infty$ with $t=O(1)$. The structure of the expansion in region $\mathbf{II_L}$, for $(-y)\gg 1$, (given by \eeref{u_2l}) suggests that in region $\mathbf{III_L}$ we write

\beq \eelab{exp_u3l}
u(y,t) = 1 - e^{-y^2 \Phi(y,t)} ,
\eeq
and expand in the form, 
\beq \eelab{exp_u3la}
\Phi(y,t) = \Phi_0(t) + \frac{1}{y} \Phi_1(t) + \frac{ \ln (- y )}{y^2} \Phi_2(t) + \frac{1}{y^2} \Phi_3(t) + o\left(y^{-2} \right),
\eeq
as $y \to - \infty$ with $t=O(1)$  and $\Phi_0(t) > 0$. On substitution from  expansions \eeref{exp_u3l} and \eeref{exp_u3la} into equation \eeref{QIVPa} we obtain a system of equations at successive orders which we solve in turn to give 
\begin{subequations}  \eelab{fun_3l}
 \begin{align}
& \Phi_0(t)=\frac{1}{(4t+C_0)}, \qquad \Phi_1(t)=\frac{(2s(t) + C_1)}{(4t + C_0)}, \qquad \Phi_2(t) = C_2,  \\
& \dot \Phi_3(t) = \dot{s}(t) \left( \frac{2s(t) + C_1}{4t + C_0} \right) + \frac{(2 + 4 C_2)}{(4t + C_0)} - \left( \frac{2s(t) + C_1}{4t + C_0} \right)^2 -f'(1), \eelab{fun_3la}
\end{align} 
\end{subequations}
where $C_0$, $C_1$,  $C_2$ and the constant associated with integrating equation \eeref{fun_3la}, $C_3$, are   constants to be determined. 
Note that $\Phi_1(t)$ and $\Phi_3(t)$ both depend on the function $s(t)$ which remains undetermined  when $t=O(1)$. We now match the expansion in region $\mathbf{III_L}$, given by substituting expressions \eeref{exp_u3la}  and \eeref{fun_3l}  into \eeref{exp_u3l} (as $t \to 0^+$), with expansion \eeref{u_2l} in region $\mathbf{II_L}$ (as $y \to -\infty$). On applying Van Dyke's matching principle \cite{VanDyke1975} we find
\beq
C_0=0, \qquad C_1=0, \qquad C_2 = -1, \qquad C_3 =\frac{1}{2} \ln \pi .
\eeq 
Thus, the expansion in region $\mathbf{III_L}$ is given by 
\beq
u(y,t) = 1  - \exp \Bigg( -\frac{y^2}{4t}  - y \frac{ s(t)}{2t}   - \ln (- y )   - \left( \frac{s(t)^2}{4t}  - \frac{1}{2} \ln t - f'(1) t + \frac{1}{2} \ln \pi  \right) +o(1) \Bigg),  \eelab{u_3l}
\eeq
as $y \to - \infty$ with $t=O(1)$. Furthermore, we note that the uniformity of expansion \eeref{u_3l}  as $y \to - \infty$ when $t \gg 1$ is dependent on the order of $s(t)$ as $t \gg 1$. This will be discussed further in section \ssref{large_time} when we investigate the asymptotic solution to QIVP as $t \to \infty$.

\subsection{Region  $\mathbf{III_R}$}
 \sslab{large_spatialb}
We next consider the corresponding region $\mathbf{III_R}$ where we determine the structure of the solution to QIVP as $y \to \infty$ with $t=O(1)$. The structure of the expansion in region $\mathbf{II_R}$, for $y \gg 1$, (given by \eeref{u_2l}) suggests that in region $\mathbf{III_R}$ we write
\beq \eelab{exp_u3r}
u(y,t) = e^{-y^2 \bar{\Phi}(y,t)},
\eeq
and expand in the form, 
\beq \eelab{exp_u3ra}
\bar{\Phi}(y,t) = \bar{\Phi}_0(t) + \frac{1}{y} \bar{\Phi}_1(t) + \frac{ \ln y}{y^2} \bar{\Phi}_2(t) + \frac{1}{y^2} \bar{\Phi}_3(t) + o\left(y^{-2} \right),
\eeq
 as $y \to \infty$ with $t=O(1)$ and $\bar{\Phi}_0(t) > 0$. On substitution from expansions \eeref{exp_u3r} and \eeref{exp_u3ra} into equation \eeref{QIVPa} we obtain a system of equations at successive orders of $y$ which we solve in turn to give 
\begin{subequations}  \eelab{fun_3r}
 \begin{align}
& \bar{\Phi}_0(t)=\frac{1}{(4t+\bar{C}_0)}, \qquad \bar{\Phi}_1(t)=\frac{(2s(t) + \bar{C}_1)}{(4t + \bar{C}_0)}, \qquad \bar{\Phi}_2(t) = \bar{C}_2,  \\
& \dot{\bar{\Phi}}_3(t) = \dot{s}(t) \left( \frac{2s(t) + \bar{C}_1}{4t + \bar{C}_0} \right) + \frac{(2 + 4 \bar{C}_2)}{(4t + \bar{C}_0)} - \left( \frac{2s(t) + \bar{C}_1}{4t + \bar{C}_0} \right)^2, \eelab{fun_3ra}
\end{align} 
\end{subequations}
 where $\bar{C}_0$, $\bar{C}_1$, $\bar{C}_2$ and the constant associated with integrating equation \eeref{fun_3ra}, $\bar{C}_3$, are   constants to be determined. %
We now match the expansion in region $\mathbf{III_R}$, given by substituting expressions \eeref{fun_3r} and \eeref{exp_u3ra} into \eeref{exp_u3r} (as $t \to 0^+$), with expansion 
 \eeref{u_2r}
 in region $\mathbf{II_R}$ (as $y \to \infty$). On applying Van Dyke's matching principle \cite{VanDyke1975} we find
\beq
\bar{C}_0=0, \qquad \bar{C}_1=0, \qquad \bar{C}_2 = -1, \qquad \bar{C}_3 =\frac{1}{2} \ln \pi .
\eeq 
Thus, the expansion in region $\mathbf{III_R}$ is given by 
\beq \eelab{u_3r}
u(y,t) =  \exp \Bigg( -\frac{y^2}{4t}  - y \frac{ s(t)}{2t}   - \ln y -  \left( \frac{s(t)^2}{4t}  - \frac{1}{2} \ln t +  \frac{1}{2} \ln \pi  \right) +o(1) \Bigg),
\eeq
as $y \to  \infty$ with $t=O(1)$. As before, the uniformity of expansion \eeref{u_3r} as 
$y \to \infty $  when $t \gg 1$ is dependent on the order of $s(t)$ as $t \gg 1$.
Finally, we are now in a position to consider the structure of the solution to QIVP 
as $t\to\infty$.

\section{Asymptotic solution to QIVP as $t \to \infty $} \sslab{large_time}
We now develop the structure of the solution to QIVP as $t \to \infty$. Guided by the numerical results in section \ssref{Numerical_results}, we anticipate that 
\beq 
s(t) = \sum^3_{i=0} c_i \phi_i(t)  + o(\phi_3(t))  \quad \mbox{as} \quad  t \to \infty, \eelab{exp_large_s}
\eeq
where $\phi_0(t)=t$, $\phi_1(t)$, $\phi_2(t)=1$ and $\phi_3(t)$ 
 are a gauge sequence as $t \to \infty$,  and the constants $c_0$, $c_1$, $c_2$,  $c_3$ are to be determined, with $c_0>0$. 
We begin by developing the structure of the solution to QIVP as $t \to \infty$ at leading order, uniform for $y \in \mathbb{R}$. We anticipate that the structure of the solution to QIVP as $t \to \infty$ will have two principal asymptotic regions in $y<0$, and two principal asymptotic regions in $y>0$. An examination of the leading order balances in the exponent of expansions \eeref{u_3l} and \eeref{u_3r} when $t \gg 1$ (using \eeref{exp_large_s}), together with the connection conditions \eeref{QIVPd} and  \eeref{QIVPe} determine the principal asymptotic structure as: 
\begin{subequations}
 \begin{align}
& \mbox{region $\mathbf{IV_L}$}:  \quad y = O(t)^-  \mbox{ with }  u=1+o(1) \mbox{  as  } t \to \infty,  \eelab{region_4l}  \\
& \mbox{region $\mathbf{IV_R}$}:  \hspace{0.31cm}  y = O(t)^+  \mbox{ with }  u=o(1) \mbox{  as  } t \to \infty, \eelab{region_4r} \\
& \mbox{region $\mathbf{V_L}$}:  \hspace{0.52cm}   y = O(1)^- \mbox{ with }  u=O(1) \mbox{  as  } t \to \infty, \eelab{region_5l}  \\
& \mbox{region $\mathbf{V_R}$}:  \hspace{0.47cm}   y = O(1)^+ \mbox{ with }  u=O(1) \mbox{  as  } t \to \infty.   \eelab{region_5r}
\end{align}  
\end{subequations}
\subsection{Regions  $\mathbf{IV_L}$,  $\mathbf{V_L}$, $\mathbf{IV_R}$
and  $\mathbf{V_R}$}

The expansion \eeref{u_3l} in region $\mathbf{III_L}$ will remain uniform for $t \gg 1$ provided that $(-y)\gg t$, but fails when $y=O(t)^-$ as $t \to \infty$. Hence, we begin in region $\mathbf{IV_L}$, in which, via \eeref{region_4l}, we introduce the scaled coordinate $w=\frac{y}{t}=O(1)^-$ as $t \to \infty$.
The structure of the expansion in region $\mathbf{III_L}$, for $t \gg 1$, (given by \eeref{u_3l}) suggests that in region $\mathbf{IV_L}$,  we write
\beq \eelab{exp_u4l}
u(w,t) = 1 - \exp \Big( -t \left( G_0(w) + o(1) \right) \Big),
\eeq
as $t \to \infty$ with $w=O(1)^-$ and $G_0(w) > 0$. On substitution of expansions  \eeref{exp_large_s} and \eeref{exp_u4l} into equation \eeref{QIVPa} we obtain the following boundary value problem, namely,
\begin{subequations} \eelab{BVP_4l}
 \begin{align}
& \left( G_0' \right)^2 - (w + c_0) G_0' + G_0  = -  f'(1) ,  \quad w < 0 , \eelab{BVP_4la}   \\
& G_0(w) > 0, \quad  w < 0 , \eelab{BVP_4lb} \\ 
& G_0(w) \sim \left(  \frac{w + c_0}{2} \right)^2 - f'(1) \quad \mbox{as} \quad w \to - \infty, \eelab{BVP_4lc} \\
& G_0(w) = O(w) \quad \mbox{as} \quad w \to 0^-. \eelab{BVP_4ld}
\end{align} 
\end{subequations}
Here condition \eeref{BVP_4lc} represents the matching condition between expansion \eeref{exp_u4l} in region $\mathbf{IV_L}$ when $(- w ) \gg 1$, and expansion \eeref{u_3l} in region $\mathbf{III_L}$ as $t \to \infty$ 
with $(-y)\gg t$
whilst condition \eeref{BVP_4ld} represents the matching condition between expansion \eeref{exp_u4l} in region $\mathbf{IV_L}$ when $w=O(t^{-1})^-$, and region $\mathbf{V_L}$ when $y=O(t)^-$ via \eeref{region_5l}. Equation \eeref{BVP_4la} has a family of linear solutions 
\beq \eelab{BVP_4l_sol1}
G_0(w)=a_1 ( w + c_0 -a_1) - f'(1) \quad \forall w < 0 ,
\eeq
for any $a_1 \in \mathbb{R}$, and an envelope solution 
\beq \eelab{BVP_4l_sol2}
G_0(w)= \left( \frac{w+c_0}{2} \right)^2 - f'(1) \quad \forall w < 0 .
\eeq
It is also possible for a combination of \eeref{BVP_4l_sol1} and \eeref{BVP_4l_sol2} to represent `envelope-linear' solutions to equation \eeref{BVP_4la}, which also remain continuous and differentiable. Applying the matching conditions \eeref{BVP_4lc} and \eeref{BVP_4ld} determines that for each $c_0 >0$, the solution to the boundary value problem \eeref{BVP_4l} is given by the `envelope-linear' solution
\beq \eelab{BVP_4l_sol3}
G_0(w) = \begin{cases}  \left( \frac{w+c_0}{2} \right)^2 - f'(1), \quad &  w < - \sqrt{c_0^2 -  4 f'(1)}  ,\\
 \left( \frac{c_0 - \sqrt{c_0^2 -  4 f'(1)}}{2} \right) w , \quad &  -  \sqrt{c_0^2 -  4 f'(1)} \leq  w < 0.  \end{cases} \\
\eeq
A sketch of $G_0(w)$, for a fixed $c_0>0$, is given in Figure \ffref{region_4}(a). For completeness we note that although $G_0(w)$ and $G_0'(w)$ are continuous, $G_0''(w)$ is discontinuous at the point $w= - \sqrt{c_0^2 -  4 f'(1)}$. Therefore, a thin transition region must exist about the point $w= - \sqrt{c_0^2 -  4 f'(1)}$ where the 
second derivative in equation \eeref{QIVPa} is retained at leading order to smooth out this discontinuity. Moreover, region $\mathbf{IV_L}$ will then be replaced by three regions, namely, region $\mathbf{IV_L^a}$, with $ - \infty < w < - \sqrt{c_0^2 -  4 f'(1)} - \textcolor{blue}{o(1)^+}$, region $\mathbf{T_L}$, a thin transition region about the point $w=-  \sqrt{c_0^2 -  4 f'(1)}$ and region $\mathbf{IV_L^b}$, with $-  \sqrt{c_0^2 -  4 f'(1)} +  \textcolor{blue}{o(1)^+} < w < 0$. As we are only interested in the leading order structure in each expansion for now, we will return to consider these regions in more detail  in \S \ref{three_L}.

   \begin{figure}[t]
     \centering
     \begin{minipage}{0.48\textwidth}
  	     \includegraphics[width=\textwidth]{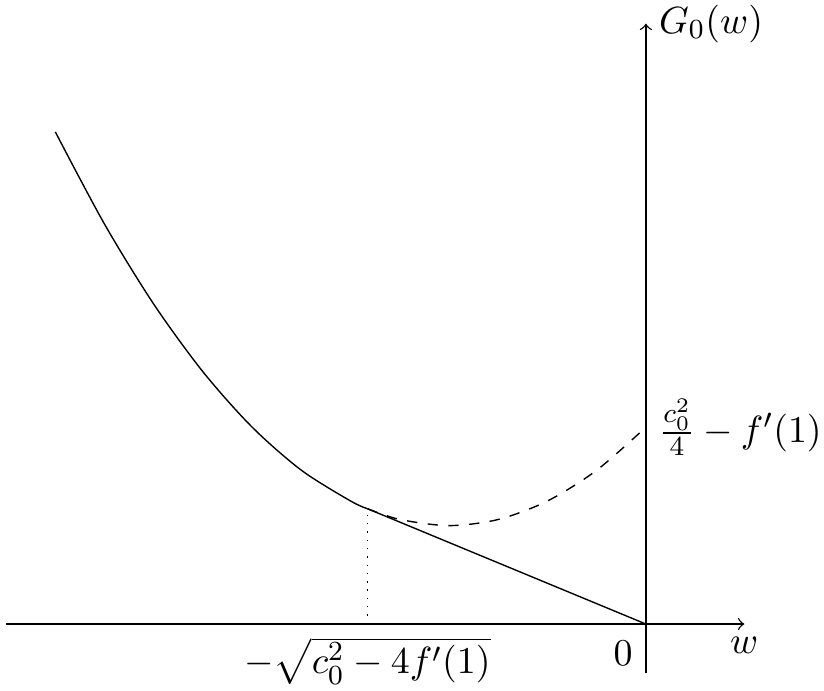}\\
  	     \centering  (a)
     \end{minipage}
     \hfill
     \begin{minipage}{0.48\textwidth}
  	   \includegraphics[width=\textwidth]{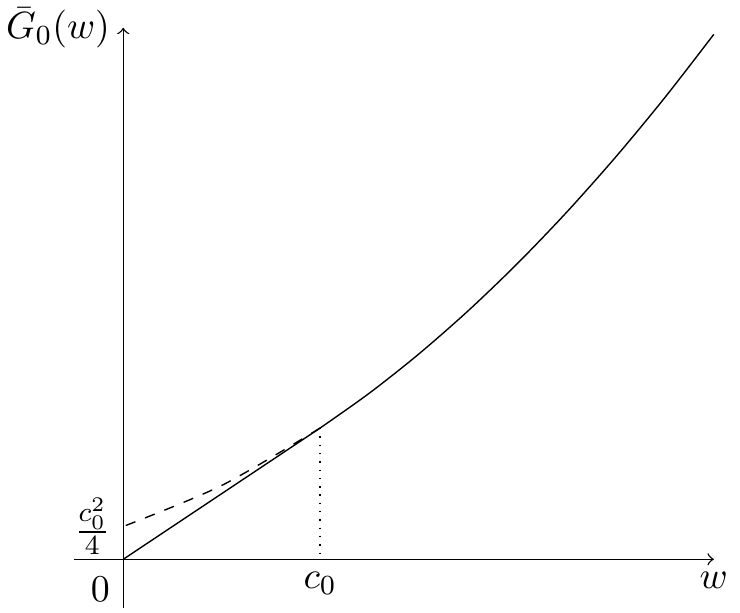}\\
    \centering  (b)
     \end{minipage}
   \caption{
    (a) A sketch of the leading order term $G_0(w)$ in the exponential expansion of the large-time solution to QIVP when $y=O(t)^-$. 
	(b) A sketch of the leading order term $\bar{G}_0(w)$ in the exponential expansion of the large-time solution to QIVP when $y=O(t)^+$.
 \fflab{region_4}
  }
   \end{figure}

Now, as $w \to 0^-$ we move out of region $\mathbf{IV_L}$ and into region $\mathbf{V_L}$, in which, via \eeref{region_5l}, $u=O(1)$ with $y=O(1)^-$ as $ t \to \infty$. In this region we therefore expand as 
\beq \eelab{exp_u5l}
u(y,t) = \hat{u}_{L0}(y) + O(\psi_L(t)) \quad \mbox{as} \quad  t \to \infty, 
\eeq
with $y =O(1)^-$, $\hat{u}_{L0}(y) > 0$  (\cite{Tisbury_etal}, equation (22b)) and where $\psi_L(t)=o(1)$ as $t \to \infty$. On substitution from expansions \eeref{exp_large_s} and \eeref{exp_u5l} into equation \eeref{QIVPa}, we obtain at leading order as $t \to \infty$, 
\begin{subequations} \eelab{BVP5l}
\beq
\hat{u}_{L0}'' + c_0 \hat{u}_{L0} ' + f(\hat{u}_{L0})=0,  \eelab{BVP5la} 
\eeq
which must be solved subject to the boundary condition \eeref{QIVPd} at $y=0$, together with the matching condition with region $\mathbf{IV_L}$ as $y \to - \infty$. Using \eeref{exp_u5l} and \eeref{BVP_4l_sol3}, these conditions require,
  \begin{align} 
& \hat{u}_{L0} (0^-) = u_c,  \\
& \hat{u}_{L0} (y) \to 1 \quad \mbox{as} \quad y \to - \infty .
\end{align} 
\end{subequations}
Due to the coupling condition \eeref{QIVPe} across $y=0$, it is necessary now to formulate the leading order problem in the corresponding regions when $y>0$ as $t \to \infty$.

The expansion \eeref{u_3r} in region $\mathbf{III_R}$ will remain uniform for $t \gg 1$ provided that $y \gg t$, but fails when $y=O(t)^+$ as $t \to \infty$. Hence, we now consider region $\mathbf{IV_R}$, in which, via \eeref{region_4r}, we introduce the scaled coordinate $w=\frac{y}{t}=O(1)^+$ as $t \to \infty$. 
The structure of the expansion in region $\mathbf{III_R}$, for $t \gg 1$, (given by \eeref{u_3r}) suggests that in region $\mathbf{IV_R}$,   we write
\beq \eelab{exp_u4r}
u(w,t) = \exp \Big( -t \left(  \bar{G}_0(w) +o(1) \right) \Big),
\eeq
as $t \to \infty$ with $w=O(1)^+$ and $\bar{G}_0(w) > 0$. On substitution of expansion \eeref{exp_u4r} into equation \eeref{QIVPa} we obtain the following boundary value problem, namely,
\begin{subequations} \eelab{BVP_4r}
  \begin{align}
& \left(\bar{G}_0' \right)^2 - (w + c_0) \bar{G}_0' + \bar{G}_0  = 0 ,  \quad w > 0 , \eelab{BVP_4ra}   \\
& \mbox{ } \bar{G}_0(w) > 0, \quad  w > 0 , \eelab{BVP_4rb} \\ 
& \mbox{ } \bar{G}_0(w) \sim \left(  \frac{w + c_0}{2} \right)^2  \quad \mbox{as} \quad w \to  \infty, \eelab{BVP_4rc}  \\
& \mbox{ } \bar{G}_0(w) = O(w)  \quad \mbox{as} \quad w \to  0^+. \eelab{BVP_4rd}
\end{align} 
\end{subequations}
Here condition \eeref{BVP_4rc} represents the matching condition between expansion \eeref{exp_u4r} in region $\mathbf{IV_R}$ when $ w \gg 1$, and expansion \eeref{u_3r} in region $\mathbf{III_R}$ as $t\to\infty$ when $y \gg t$ whilst condition \eeref{BVP_4rd} represents the matching condition between expansion \eeref{exp_u4r} in region $\mathbf{IV_R}$ when $ w=O(t^{-1})^+$, and region $\mathbf{V_R}$ when $y=O(t)^+$ via \eeref{region_5r}. For each $c_0> 0$, the boundary value problem \eeref{BVP_4r} has the unique solution 
\beq \eelab{BVP_4r_sol}
\bar{G}_0(w) = \begin{cases} \left( \frac{w+c_0}{2} \right)^2 , \quad &  w > c_0 ,\\
  \; c_0 w , \quad &  0 < w \leq c_0.  \end{cases} \\
\eeq
A sketch of $\bar{G}_0(w)$ for a fixed $c_0>0$ is given in Figure \ffref{region_4}(b). For completeness we note that although $\bar{G}_0(w)$ and $\bar{G}_0'(w)$ are continuous, $\bar{G}_0''(w)$ is discontinuous at the point $w=c_0$. Hence, a thin transition region about the point $w=c_0$ is required in which the second derivative in equation \eeref{QIVPa} is retained at leading order to smooth out the discontinuity. This requires that region $\mathbf{IV_R}$ is replaced by three regions, namely, region $\mathbf{IV_R^a}$, with $c_0 + o(1) < w < \infty$, region $\mathbf{T_R}$, a thin transition region about the point $w=c_0$ and region $\mathbf{IV_R^b}$, with $0 < w < c_0 -o(1)$. As before, we will consider these regions in more detail 
in \S \ref{three_R}.

Now, as $w \to 0^+$ we move out of region $\mathbf{IV_R}$ and into region $\mathbf{V_R}$, in which, via \eeref{region_5r}, $u=O(1)$ and $y=O(1)^+$ as $ t \to \infty$. In this region we must therefore expand as 
\beq \eelab{exp_u5r}
u(y,t) = \hat{u}_{R0}(y) + O(\psi_R(t)) \quad \mbox{as} \quad  t \to \infty, 
\eeq
with $y =O(1)^+$, $\hat{u}_{R0}(y) > 0$ (\cite{Tisbury_etal}, equation (22b)) and $\psi_R(t)=o(1)$ as $t \to \infty$. On substitution from expansions \eeref{exp_large_s} and \eeref{exp_u5r} into equation \eeref{QIVPa}, we obtain at leading order as $t \to \infty$, 
\begin{subequations} \eelab{BVP5r}
\beq
\hat{u}_{R0}'' + c_0 \hat{u}_{R0} ' = 0, \eelab{BVP5ra} 
\eeq
which must be solved subject to the boundary condition \eeref{QIVPd} at $y=0$, together with the matching condition with region $\mathbf{IV_R}$ as $y \to  \infty$. Using \eeref{exp_u5l} and \eeref{BVP_4l_sol3}, these conditions require,
  \begin{align} 
& \hat{u}_{R0} (0^+)  = u_c,  \\
& \hat{u}_{R0} (y)  \to 0 \quad \mbox{as} \quad y \to  \infty .
\end{align}  
\end{subequations}
Finally, the boundary value problems \eeref{BVP5l} and \eeref{BVP5r} must be solved subject to the coupling condition \eeref{QIVPe} across $y=0$, which requires
\beq \eelab{BVP5c}
\hat{u}_{L0}' (0^-) = \hat{u}_{R0}' (0^+).
\eeq

The coupled nonlinear boundary value problem, given by \eeref{BVP5l}, \eeref{BVP5r} and \eeref{BVP5c}, across regions $\mathbf{V_L}$ and $\mathbf{V_R}$ is precisely the nonlinear boundary value problem satisfied by the PTW structure considered in Part I with $v$ replaced by $c_0$. Thus, we immediately conclude that
\begin{subequations} \eelab{BVP_TW}
  \begin{align}
& \hat{u}_{R0}(y)=U_T(y), \quad y \geq 0,   \eelab{BVP_TWa}   \\
& \hat{u}_{L0}(y)=U_T(y), \quad y < 0,   \eelab{BVP_TWb} 
\end{align}  
and that $c_0$ is now determined as,
\beq \eelab{c0}
c_0 = v^*(u_c),
\eeq
\end{subequations}
where $U_T:\mathbb{R} \to \mathbb{R}$ is the PTW solution to QIVP at cut-off $u_c \in (0,1)$, which has propagation speed $v^*(u_c)$. 
For convenience, we recall from Theorem 1.1 of Part I that 
 \begin{subequations}\eelab{exp_TW}   
 \begin{align}
 & U_T(y)= u_c e^{- v^*(u_c)y}  \quad \forall y \in [0, \infty),  \eelab{exp_TWa}  \\
 \intertext{and}
 & U_T(y) \sim 1 - A_{-\infty} e^{\lambda_+ ( v^*(u_c))y} \quad \mbox{as} \quad y \to - \infty,   \eelab{exp_TWb}   
\end{align} 
\end{subequations}
where $\lambda_+(v^*(u_c)) =  \frac{1}{2} \left( - v^*(u_c) + \sqrt{v^*(u_c)^2 - 4 f'(1) } \right)$, and $A_{-\infty}$ is a global constant depending upon $u_c$. 
 This completes the asymptotic structure of the solution to QIVP as $t \to \infty$ at leading order.

\subsection{Regions $\mathbf{IV_R^a}$,
$\mathbf{T_R}$, 
 $\mathbf{IV_R^b}$ and $\mathbf{V_R}$}\label{three_R}

To develop the solution to QIVP to higher order we must first return to region $\mathbf{T_R}$, the localised transition region in which $w=v^*(u_c)+o(1)$ as $t \to \infty$. It follows from the leading order term in the expansion in region $\mathbf{IV_R}$ (given by \eeref{BVP_4r_sol}, \eeref{BVP5r} and \eeref{c0}) that to examine region $\mathbf{T_R}$ we must introduce the scaled coordinate $\zeta=(w-v^*(u_c))t^{\frac{1}{2}}$ 
and expand $u(\zeta,t)$ in the form
\beq \eelab{exp_TR}
u(\zeta,t) = \left( \bar{F}_0(\zeta) + o(1) \right) \exp \bigg(-t v^*(u_c)^2 - t^{\frac{1}{2}} \zeta v^*(u_c) \bigg),
\eeq
as $t \to \infty$ with $\zeta=O(1)$ and $\bar{F}_0(\zeta)>0$. On substitution of expansions \eeref{exp_TR} and \eeref{exp_large_s} into equation \eeref{QIVPa} we obtain  
\beq
t \dot\phi_1(t) \left( v^*(u_c) c_1 \bar{F}_0 \right) + \left( - \frac{1}{2} \zeta \bar{F}_0' - \bar{F}_0'' \right) + o(1)= 0,  \qquad - \infty < \zeta < \infty.
\eeq
The only non-trivial dominant balance requires that we set, without loss of generality
\beq \eelab{phi1}
\phi_1(t) = \ln t.
\eeq
Thus, the leading order equation in region $\mathbf{T_R}$ is given by 
\beq \eelab{TR1}
\bar{F}_0'' + \frac{1}{2} \zeta \bar{F}_0' - \gamma \bar{F}_0= 0, \qquad - \infty < \zeta < \infty,
\eeq
with $\gamma = v^*(u_c) c_1$. To obtain the full boundary value problem for $\bar{F}_0(\zeta)$ we require matching conditions as $\zeta \to - \infty$ with region $\mathbf{IV_R^b}$ and as $\zeta \to  \infty$ with region $\mathbf{IV_R^a}$. Therefore, we next return to region $\mathbf{IV_R^b}$. The structure of the expansion in region $\mathbf{V_R}$, for $y \gg 1$, (given by \eeref{exp_u5r}, \eeref{BVP_TWa} and \eeref{exp_TWa}) dictates that in region $\mathbf{IV_R^b}$ we expand in the form
\beq \eelab{exp_u4r_full}
u(w,t) = \exp \left( -t \left( v^*(u_c) w - \frac{1}{t} \hat{G}(w) + o \left( \frac{1}{t} \right) \right) \right),
\eeq
as $t \to \infty$ with $O(t^{-1}) < w < v^*(u_c) - O(t^{-\frac{1}{2}})$. We substitute expansion \eeref{exp_u4r_full} into equation \eeref{QIVPa} to obtain, on solving at each order in turn, 
\beq \eelab{exp_u4r_full2}
u(w,t) = \exp \bigg( -t  v^*(u_c) w +  v^*(u_c) c_1 \ln \big( v^*(u_c)-w \big)+ \bar{d} + o \left(1 \right)  \bigg),
\eeq
as $t \to \infty$ with $O(t^{-1}) < w < v^*(u_c) - O(t^{-\frac{1}{2}})$ and where the constants $c_1$ and $\bar{d}$ are to be determined. On matching expansion \eeref{exp_u4r_full2} in region $\mathbf{IV_R^b}$ (as $w \to v^*(u_c)^-$) with expansion \eeref{exp_TR} in region $\mathbf{T_R}$ (as $\zeta \to - \infty$), via Van Dyke's matching principle \cite{VanDyke1975}, we readily obtain that 
\beq \eelab{c1}
c_1 =0,
\eeq 
after which we must have
\beq \eelab{TR2}
\bar{F}_0(\zeta) = e^{\bar{d}} + o(1) \quad \mbox{as} \quad \zeta \to - \infty.
\eeq
To determine $\bar{d}$ we next match expansion \eeref{exp_u4r_full2} (with \eeref{c1}) in region $\mathbf{IV_R^b}$ (as $w \to 0^+$) with expansion \eeref{exp_TWa} in region $\mathbf{V_R}$ (as $y \to \infty$). On applying Van Dyke's matching principle \cite{VanDyke1975}, we require that 
\beq \eelab{bar_d}
\bar{d} = \ln u_c .
\eeq 
Thus, via \eeref{exp_u4r_full2}, \eeref{c1} and \eeref{bar_d}, the expansion in region $\mathbf{IV_R^b}$ is given by
\beq \eelab{u4r_full}
u(w,t) =  \exp \bigg( -t  v^*(u_c) w + \ln u_c + o(1)  \bigg),
\eeq
as $t \to \infty$ with $O(t^{-1}) < w < v^*(u_c) - O(t^{-\frac{1}{2}})$. In addition \eeref{TR2} becomes 
\beq \eelab{TR3}
\bar{F}_0(\zeta) = u_c + o(1) \quad \mbox{as} \quad \zeta \to - \infty.
\eeq
We next consider region $\mathbf{IV_R^a}$. The structure of the expansion in region $\mathbf{III_R}$, as $t \to \infty$ with $y\gg t$, (given by \eeref{u_3r}) and the form of $s(t)$ as $t \to \infty$ (given by \eeref{exp_large_s} with $c_1$ now determined by \eeref{c1}), suggests that in region $\mathbf{IV_R^a}$ we write
\beq \eelab{exp_u_4ra}
u(w,t) = e^{-t \bar{G}(w,t)},
\eeq
and expand in the form,
\beq \eelab{exp_u_4ra2}
\bar{G}(w,t) = \left( \frac{w+v^*(u_c)}{2} \right)^2 + \frac{\ln t }{t}\bar{G}_1(w) + \frac{1}{t} \bar{G}_2(w) + o(t^{-1}),
\eeq
as $t \to \infty$ with $w > v^*(u_c)+ O(t^{-\frac{1}{2}})$. On substitution from  \eeref{exp_u_4ra} and \eeref{exp_u_4ra2} into equation \eeref{QIVPa} we obtain a series of boundary value problems which we solve at each order of $t$ in turn to obtain
\beq \eelab{exp_u_4ra3}
u(w,t) =  \exp \bigg( - t \left( \frac{w+v^*(u_c)}{2} \right)^2  - \frac{1}{2} \ln t  -  \bar{G}_2(w)+o(1) \bigg),
\eeq
as $t \to \infty$ with $w > v^*(u_c)+ O(t^{-\frac{1}{2}})$ and where the function $\bar{G}_2(w)$ is indeterminate, being globally dependent on the evolution at earlier stages when $t=O(1)$ and $y=O(1)$. However, to match with expansion $\mathbf{III_R}$ 
(as $t\to\infty$ with $y\gg t$), we require
\beq \eelab{match_con_4r1}
\bar{G}_2(w) \sim c_2 \left( \frac{w + v^*(u_c) }{2} \right) + \ln w + \frac{1}{2} \ln \pi \quad \mbox{as} \quad w \to \infty.
\eeq
In addition the structure of the expansion in region $\mathbf{T_R}$, as given by \eeref{exp_TR}, requires, for matching to be possible, that,
\beq \eelab{match_con_4r2}
\bar{G}_2(w) \sim \bar{\alpha}_1 \ln \big( w - v^*(u_c) \big) + \bar{\alpha}_2  \quad \mbox{as} \quad w \to v^*(u_c)^+,
\eeq
for some constants $\bar{\alpha}_1, \bar{\alpha}_2$ to be determined. We now match in detail the expansion in region $\mathbf{IV_R^a}$, given by  \eeref{exp_u_4ra3} and \eeref{match_con_4r2} (as $w \to v^*(u_c)^+$), with expansion \eeref{exp_TR} in region $\mathbf{T_R}$ (as $\zeta \to \infty$). On applying Van Dyke's matching principle \cite{VanDyke1975} we find that
\beq
\bar{\alpha}_1 = 1, 
\eeq
after which,
\beq \eelab{TR4}
\bar{F}_0(\zeta) = \bar{\sigma} \zeta^{-1} e^{- \frac{\zeta^2}{4}} \left(1 + o(1) \right) \quad \mbox{as} \quad \zeta \to  \infty,
\eeq
where $\bar{\sigma} =e^{-\bar{\alpha}_2}$. Hence, on collecting \eeref{TR1}, \eeref{c1}, \eeref{TR3} and \eeref{TR4} we obtain the boundary value problem in region $\mathbf{T_R}$ for $\bar{F}_0(\zeta)$ as,
\begin{subequations} \eelab{BVP_TR}
 \begin{align}
& \bar{F}_0'' + \frac{1}{2} \zeta \bar{F}_0' = 0, \qquad - \infty < \zeta < \infty ,   \eelab{BVP_TRa}    \\
& \bar{F}_0(\zeta)> 0, \qquad - \infty < \zeta < \infty  ,  \\
& \bar{F}_0(\zeta) = \bar{\sigma} \zeta^{-1} e^{- \frac{\zeta^2}{4}} \left(1 + o(1) \right) \quad \mbox{as} \quad \zeta \to  \infty , \eelab{BVP_TRc}     \\
& \bar{F}_0(\zeta) = u_c + o(1) \quad \mbox{as} \quad \zeta \to  -\infty . \eelab{BVP_TRd}
\end{align} 
\end{subequations}
This boundary value problem has a solution only when 
\beq \eelab{bar_sigma}
\bar{\sigma}= \frac{u_c}{\sqrt{\pi}},
\eeq
with the solution being unique, and given by, 
\beq
\bar{F}_0(\zeta) = \frac{1}{2} u_c \erfc \left( \frac{\zeta}{2} \right)
\quad \forall - \infty < \zeta < \infty. 
\eeq
It follows from \eeref{bar_sigma} that
\beq
\bar{\alpha}_2 = - \ln \frac{u_c}{\sqrt{\pi}}.
\eeq

\begin{figure}
\centering
{\includegraphics[width=0.95\linewidth]{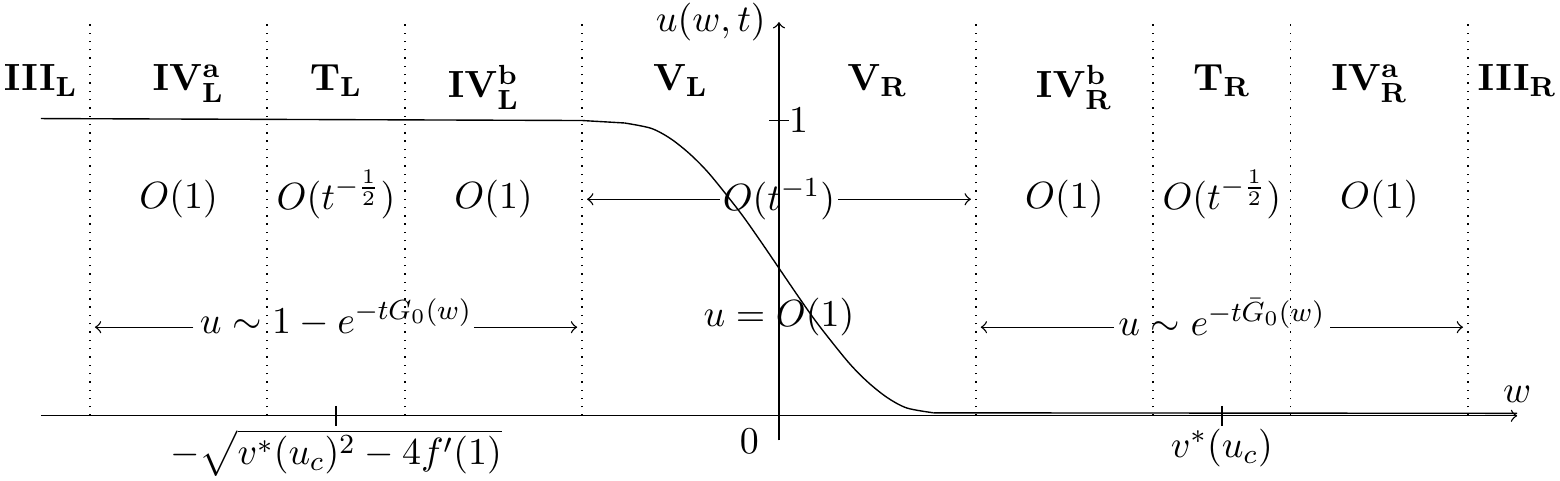}} 
\caption{A schematic representation of the location and thickness of the asymptotic regions in the solution to QIVP as $t \to \infty$. Here the the leading order terms in the exponential form of the solution $G_0(w)$ and $\bar{G}_0(w)$ are given by \eeref{BVP_4l_sol3} and \eeref{BVP_4r_sol}, respectively. Additionally, there are thin transition regions at $w= - \sqrt{ v^*(u_c)^2 - 4f'(1)}$ and at $w=v^*(u_c)$. Note that regions $\mathbf{III_L}$ and $\mathbf{III_R}$ are far field regions for $\lvert w \rvert \gg1$ as $t \to \infty$.}
\fflab{large_time}
\end{figure}

It is now instructive to summarize the structure in regions $\mathbf{IV_R^a}, \mathbf{T_R}$ and $\mathbf{IV_R^b}$. The expansion in region $\mathbf{IV_R^a}$ is given by \eeref{exp_u_4ra3} together with the asymptotic conditions 
 \beq \eelab{match_con_4l1}
\bar{G}_2(w) \sim  \begin{cases}  \ln \big( w - v^*(u_c)  \big) - \ln \frac{u_c}{\sqrt{\pi}} , \quad &  \mbox{as} \quad w \to v^*(u_c)^+,\\
  c_2\left( \frac{w + v^*(u_c) }{2} \right) + \ln w + \frac{1}{2} \ln \pi , \quad & \mbox{as} \quad w \to \infty, \end{cases} \\ 
\eeq
whilst in region $\mathbf{T_R}$
\beq \eelab{u_TR}
u(\zeta,t) = \left( \frac{1}{2} u_c \erfc \left( \frac{\zeta}{2} \right) + o(1) \right) \exp \bigg( -t v^*(u_c)^2 - t^{\frac{1}{2}} \zeta v^*(u_c) \bigg) ,
\eeq
as $t \to \infty$ with $\zeta = O(1)$, and in region $\mathbf{IV_R^b}$
\beq \eelab{u_4rb}
u(w,t) = \exp \bigg( -t v^*(u_c) w + \ln u_c + o(1) \bigg) ,
\eeq
as $t \to \infty$ with $O(t^{-1}) < w < v^*(u_c) - O(t^{-\frac{1}{2}})$. A schematic representation of the location and thickness of the asymptotic regions as $t \to \infty$ is given in Figure \ffref{large_time}.

We next consider the structure of the expansion in region $\mathbf{T_R}$ in more detail. Via \eeref{u_TR}, we observe that for $(- \zeta) \gg 1$,
\beq \eelab{exp_uTR_match}
u(\zeta,t) \sim \exp \left(-t v^*(u_c)^2 - t^{\frac{1}{2}} v^*(u_c) \zeta + 
\ln \left(u_c\left(1 + \frac{1}{\sqrt{\pi}} \frac{1}{\zeta} e^{- \frac{\zeta^2}{4}} \right)\right)\right),
\eeq
as $t \to \infty$, which demands that in region $\mathbf{IV_R^b}$, to continue the expansion in \eeref{u_4rb}, we must write
\beq \eelab{exp_u_4rb}
u(w, t) = u_c e^{-t w v^*(u_c)} + t^{-\frac{1}{2}} \bar{G}(w,t) \exp \left( -\frac{t (w+v^*(u_c))^2}{4} \right),
\eeq
as $t \to \infty$ with $O(t^{-1}) < w < v^*(u_c) - O(t^{-\frac{1}{2}})$ and $\bar{G}(w,t) =O(1)$ as $t \to \infty$. On substituting from expansion \eeref{exp_u_4rb} into equation \eeref{QIVPa}, and simplifying, we obtain
\beq
\bar{G}_t -\frac{1}{2}t^{-1}\bar{G}- t^{-2} \bar{G}_{ww} =  O \left( t^{\frac{1}{2}} \dot\phi_3(t) \exp  \left( -t \left( w v^*(u_c) - \frac{(w+ v^*(u_c))^2}{4} \right) \right) \right), \eelab{4rb_1}
\eeq
as $t \to \infty$ with $O(t^{-1}) < w < v^*(u_c) - O(t^{-\frac{1}{2}})$. We will later verify that the right-hand side of equation \eeref{4rb_1} is exponentially small as $t \to \infty$ in this region. Hence, to obtain a structured balance in \eeref{4rb_1}, we must expand $\bar{G}(w,t)$ in the form
\beq \eelab{4rb_2}
\bar{G}(w,t) = \bar{G}_0(w) + t^{-1} \bar{G}_1(w) + o \left( t^{- 1} \right),
\eeq
as $t \to \infty$ with $O(t^{-1}) < w < v^*(u_c) - O(t^{-\frac{1}{2}})$ and on substitution into \eeref{4rb_1} we obtain at leading order 
\beq \eelab{4rb_3}
\bar{G}_0'' +  \bar{G}_1 =0,
\eeq
with $O(t^{-1}) < w < v^*(u_c) - O(t^{-\frac{1}{2}})$. We conclude that $\bar{G}_0(w)$ is indeterminate and represents a further globally determined function. Therefore, the expansion in region $\mathbf{IV_R^b}$ is, from equations \eeref{exp_u_4rb} and \eeref{4rb_2},
 \begin{equation}
u(w, t) =   u_c e^{-t w v^*(u_c)} + t^{-\frac{1}{2}} \bar{G}_0(w)( 1 + O(t^{-1}))  \exp \left( -\frac{t (w+v^*(u_c))^2}{4} \right),  \eelab{exp_u_4rb2}
\end{equation}
as $t \to \infty$ with $O(t^{-1}) < w < v^*(u_c) - O(t^{-\frac{1}{2}})$. We now match the expansion \eeref{exp_u_4rb2} in region $\mathbf{IV_R^b}$ (as $w \to v^*(u_c)^-$), with expansion \eeref{exp_uTR_match} in region $\mathbf{T_R}$ (as $\zeta \to - \infty$), in detail. On applying Van Dyke's matching principle \cite{VanDyke1975} we require 
\beq
\bar{G}_0(w) = - \frac{u_c}{\sqrt{\pi}} (w - v^*(u_c))^{-1} + o(w - v^*(u_c))^{-1} \quad \mbox{as} \quad w \to v^*(u_c)^-.
\eeq
We next return to region $\mathbf{V_R}$. First, a balance between expansion 
\eeref{exp_u5l} in region $\mathbf{V_L}$ and expansion \eeref{exp_u5r} in region $\mathbf{V_R}$, across the connection at $y=0$, requires 
\beq \eelab{region5_bal}
\psi_L(t)=\psi_R(t)=\psi(t),
\eeq
where $\psi(t)=o(1)$ as $t \to \infty$. Now, the induced correction term in expansion \eeref{exp_u5r} in region $\mathbf{V_R}$ from region $\mathbf{IV_R^b}$ when $0<w \ll 1$, must have, via \eeref{exp_u_4rb2},
\beq \eelab{region5_cor1}
\psi(t) = O \left( t^{\gamma} e^{- \frac{v^*(u_c)^2 t }{4}}  \right),
\eeq
as $t \to \infty$, with constant $\gamma$ to be determined. Thus, without loss of generality we set 
\beq \eelab{region5_cor2}
\psi(t) = t^{\gamma} e^{- \frac{v^*(u_c)^2 t }{4}} .
\eeq
Hence, in region $\mathbf{V_R}$ we develop expansion \eeref{exp_u5r} in the form  
\beq \eelab{exp_u5r_2}
u(y,t) = U_T(y) + t^{\gamma} e^{- \frac{v^*(u_c)^2 t }{4}} u_1(y) (1 + o(1)),
\eeq
as $t \to \infty$ with $y=O(1)^+$. On substitution of expansion \eeref{exp_u5r_2} into equation \eeref{QIVPa}, and cancelling at leading order, we obtain
\beq
 - \frac{1}{4} v^*(u_c)^2 u_1 - v^*(u_c) u_1' - u_1'' +o(1) =
   c_3 U_T'(y) t^{- \gamma} \dot\phi_3(t) e^{\frac{v^*(u_c)^2 t }{4}} , \eelab{region5_1}
\eeq
as $t \to \infty$ with $y=O(1)^+$. The non-trivial balance in \eeref{region5_1} requires that we set, without loss of generality
\beq \eelab{phi3}
\dot\phi_3(t) = t^{\gamma} e^{- \frac{v^*(u_c)^2 t }{4}},
\eeq
and we note that this now confirms that the right-hand side of \eeref{4rb_1} is exponentially small as $t \to \infty$. The corresponding problem for $u_1(y)$ is then
\begin{subequations} \eelab{BVP2_5r}
  \begin{align}
& u_1'' + v^*(u_c) u_1' + \frac{1}{4} v^*(u_c)^2  u_1 = -c_3 U_T'(y), \quad y > 0 ,   \eelab{BVP2_5ra} \\
& u_1(0^+)=0, \eelab{BVP2_5rb}
\end{align}  
\end{subequations}
where the condition \eeref{BVP2_5rb} is required for the boundary condition \eeref{QIVPd} to be satisfied. The problem for $u_1(y)$, given by \eeref{BVP2_5r}, must be solved subject to the matching condition with region $\mathbf{IV_R^b}$. Before formulating this matching condition, we consider the corresponding structure in regions $\mathbf{IV_L^a}, \mathbf{T_L}, \mathbf{IV_L^b}$ and $\mathbf{V_L}$. Thus, we now move   to region $\mathbf{IV_L^a}$.

\subsection{Regions $\mathbf{IV_L^a}$,
$\mathbf{T_L}$, 
 $\mathbf{IV_L^b}$  and $\mathbf{V_L}$}\label{three_L}

The structure of the expansion in region $\mathbf{III_L}$ as $t\to\infty$ with $(-y)\gg t$ (given by \eeref{u_3l}), the structure of $s(t)$ as $t \to \infty$ (given by  \eeref{exp_large_s} with $c_0$ and $c_1$ given by \eeref{c0} and \eeref{c1} respectively) and the leading order behaviour in regions $\mathbf{IV_L^a}$ and $\mathbf{IV_L^b}$ (given by \eeref{exp_u4l} and \eeref{BVP_4l_sol3}), suggests that in region $\mathbf{IV_L^a}$ we write 
\beq \eelab{exp_u4la}
u(w,t) = 1 - e^{-t G(w,t)},
\eeq
and expand in the form,
\beq \eelab{exp_u4la_2}
G(w,t) = \left( \frac{w+v^*(u_c)}{2} \right)^2 - f'(1) + \frac{\ln t }{t}G_1(w) + \frac{1}{t} G_2(w) + o(t^{-1}),
\eeq
as $t \to \infty$ with $  w < - \sqrt{v^*(u_c)^2 -  4 f'(1)} - O(t^{-\frac{1}{2}})$. On substitution of \eeref{exp_u4la} and expansion \eeref{exp_u4la_2} into equation \eeref{QIVPa} we obtain a sequence of boundary value problems which we solve at each order   to obtain
\beq \eelab{u4la}
u(w,t) =  1- \exp \Bigg( - t \left( \left( \frac{w+v^*(u_c)}{2} \right)^2 -f'(1) \right) - \frac{1}{2}  \ln t  -  G_2(w)+o(1) \Bigg), 
\eeq
as $t \to \infty$ with $  w < - \sqrt{v^*(u_c)^2 -  4 f'(1)} - O(t^{-\frac{1}{2}})$, and where the function $G_2(w)$ is indeterminate, being globally dependent on the evolution at earlier stages when $t=O(1)$ and $y=O(1)$. However, to match with expansion $\mathbf{III_L}$ (as $t\to\infty$ with $(-y)\gg t$), we require
\beq \eelab{match_con_4la}
G_2(w) \sim c_2 \left( \frac{w + v^*(u_c) }{2} \right) + \ln (-w)+ \frac{1}{2} \ln \pi \quad \mbox{as} \quad w \to - \infty.
\eeq
We next examine region $\mathbf{T_L}$. It follows from the structure of the expansion in region $\mathbf{IV_L^a}$, as $w \to (-  \sqrt{v^*(u_c)^2 -  4 f'(1)})^-$ (given by \eeref{u4la}), that in region $\mathbf{T_L}$ we must introduce the scaled coordinate $\zeta=\left(w+\sqrt{v^*(u_c)^2 - 4f'(1)}\right)t^{\frac{1}{2}}$ %
and expand $u(\zeta,t)$ in the form
  \begin{align} \eelab{exp_TL}
u(\zeta,t) =  & 1 - \left( F_0(\zeta) + o(1) \right) \exp \Bigg( -t \Bigg( \left( \frac{ v^*(u_c) - \sqrt{v^*(u_c)^2 - 4f'(1)}}{2} \right)^2  \nonumber  \\ 
&  -f'(1) \Bigg)   - t^{\frac{1}{2}} \zeta \left( \frac{v^*(u_c) - \sqrt{v^*(u_c)^2 - 4f'(1)}}{2} \right)  \Bigg)  ,
\end{align} 
as $t \to \infty$ with $\zeta=O(1)$. On substitution of expansion \eeref{exp_TL} into equation \eeref{QIVPa} we obtain at leading order 
\beq \eelab{ODE_TL}
F_0'' + \frac{1}{2}\zeta F_0'= 0, \qquad - \infty < \zeta < \infty.
\eeq
To obtain the full boundary value problem for $F_0(\zeta)$ we require matching conditions as $\zeta \to \pm \infty$. To that end, the structure of the expansion in region $\mathbf{T_L}$, as given by \eeref{exp_TL}, requires, for matching to be possible, with expansions \eeref{u4la}  and 
\eeref{match_con_4la} 
in region  $\mathbf{IV_L^a}$, that
\beq \eelab{match_con_TL1}
G_2(w) \sim \alpha_1 \ln \left|  w+\sqrt{v^*(u_c)^2 - 4f'(1)} \right|+ \alpha_2, 
\eeq
as $w \to (-\sqrt{v^*(u_c)^2 - 4f'(1)})^- $ for some   constants $\alpha_1, \alpha_2$ to be determined. We now match in detail the expansion in region $\mathbf{IV_L^a}$, given by \eeref{u4la} and \eeref{match_con_TL1}, as $w \to (-  \sqrt{v^*(u_c)^2 -  4 f'(1)})^-$, with expansion \eeref{exp_TL} in region $\mathbf{T_L}$, as $\zeta \to - \infty$. On applying Van Dyke's matching principle \cite{VanDyke1975} it immediately follows that  
\beq
\alpha_1 = 1,
\eeq
after which we must have 
\beq \eelab{match_con_TL2}
F_0(\zeta) = \sigma \zeta^{-1} e^{- \frac{\zeta^2}{4}} \left(1 + o(1) \right) \quad \mbox{as} \quad \zeta \to - \infty,
\eeq
where $\sigma=e^{-\alpha_2}$. We next consider the matching condition as $\zeta \to \infty$. The structure of the expansion in region $\mathbf{V_L}$, for $(-y) \gg 1$, (given by \eeref{exp_u5l}, \eeref{BVP_TWb}  and \eeref{exp_TWb}) dictates that in region $\mathbf{IV_L^b}$ we must expand in the form
\beq
u(w,t) =  1- \exp \left( -t  \left( \frac{v^*(u_c) - \sqrt{v^*(u_c)^2 - 4f'(1)}}{2} \right) w   +   \tilde{G}(w) + o \left( 1 \right) \right) , \eelab{u4lb_1}
\eeq
as $t \to \infty$ with $-  \sqrt{v^*(u_c)^2 -  4 f'(1)} + O(t^{-\frac{1}{2}}) < w < O(t^{-1})^-$. We substitute expansion \eeref{u4lb_1} into equation \eeref{QIVPa} to obtain, on solving at each order in turn, 
\beq
u(w,t) =  1- \exp \left( -t  \left( \frac{v^*(u_c) - \sqrt{v^*(u_c)^2 - 4f'(1)}}{2} \right) w  + d + o \left(1 \right)  \right), \eelab{u4lb_2}
\eeq
as $t \to \infty$ with  $-  \sqrt{v^*(u_c)^2 -  4 f'(1)} + O(t^{-\frac{1}{2}}) < w < O(t^{-1})^-$ and where the constant $d$ is to be determined. On matching expansion \eeref{u4lb_2} in region $\mathbf{IV_L^b}$ (as $w \to 0^-$) with expansion \eeref{exp_TWb} in region $\mathbf{V_L}$ (as $y \to - \infty$), via Van Dyke's matching principle \cite{VanDyke1975}, we readily obtain that 
\beq \eelab{d}
d=\ln A_{-\infty}.
\eeq
 Thus, via \eeref{u4lb_2} and \eeref{d}, the expansion in region $\mathbf{IV_L^b}$ is given by
\beq \eelab{u4lb_3}
u(w,t) = 1 - \exp \left( -t  \left( \frac{v^*(u_c) - \sqrt{v^*(u_c)^2 - 4f'(1)}}{2} \right) w  + \ln A_{-\infty} + o \left(1 \right)  \right),
\eeq
as $t \to \infty$ with $-  \sqrt{v^*(u_c)^2 -  4 f'(1)} + O(t^{-\frac{1}{2}}) < w < O(t^{-1})^-$. On matching expansion \eeref{u4lb_3} in region $\mathbf{IV_L^b}$ (as $w \to (-  \sqrt{v^*(u_c)^2 -  4 f'(1)})^-$) with expansion \eeref{exp_TL} in region $\mathbf{T_L}$ (as $\zeta \to \infty$), we obtain the condition
\beq \eelab{match_con_TL3}
F_0(\zeta) = A_{-\infty} + o(1) \quad \mbox{as} \quad \zeta \to  \infty.
\eeq
Hence, on collecting \eeref{ODE_TL}, \eeref{match_con_TL2} and \eeref{match_con_TL3} we obtain the boundary value problem in region $\mathbf{T_L}$ for $F_0(\zeta)$ as,
\begin{subequations} \eelab{BVP_TL}
  \begin{align}
& F_0'' + \frac{1}{2}\zeta F_0'= 0, \qquad - \infty < \zeta < \infty ,   \eelab{BVP_TL1}    \\
& F_0(\zeta)> 0, \qquad - \infty < \zeta < \infty,    \\
& F_0(\zeta) = \sigma \zeta^{-1} e^{- \frac{\zeta^2}{4}} \left(1 + o(1) \right) \quad \mbox{as} \quad \zeta \to  - \infty , \eelab{BVP_TL2}     \\
& F_0(\zeta) = A_{-\infty} + o(1) \quad \mbox{as} \quad \zeta \to  \infty . \eelab{BVP_TL3}
\end{align}  
\end{subequations}
This boundary value problem has a solution only when 
\beq \eelab{sigma}
\sigma= \frac{A_{-\infty}}{\sqrt{\pi}},
\eeq
with the solution being unique, and given by, 
\beq
F_0(\zeta) = \frac{1}{2} A_{-\infty} \left(1 + \erf \left( \frac{\zeta}{2} \right) \right)\qquad \forall - \infty < \zeta < \infty.
\eeq
It follows from \eeref{sigma} that
\beq
\alpha_2 = - \ln \frac{A_{-\infty}}{\sqrt{\pi}}.
\eeq

It is again instructive to summarize the structure in regions $\mathbf{IV_L^a}, \mathbf{T_L}$ and $\mathbf{IV_L^b}$. The expansion in region $\mathbf{IV_L^a}$ is given by \eeref{u4la} together with the asymptotic conditions
\begin{equation}  \eelab{match_con_4la_full}
G_2(w) \sim  \begin{cases}  
\ln \lvert w + \sqrt{v^*(u_c)^2 - 4f'(1)} \rvert - \ln \frac{A_{-\infty}}{\sqrt{\pi}}, &  \mbox{as} \; w \to \left( -\sqrt{v^*(u_c)^2 - 4f'(1)} \right)^-,\\
 \left( \frac{w + v^*(u_c) }{2} \right) + \ln \lvert w \rvert + \frac{1}{2} \ln \pi, & \mbox{as} \; w \to -\infty, \end{cases}  
\end{equation}
whilst in region $\mathbf{T_L}$,
  \begin{align} \eelab{exp_uTL}
u(\zeta,t) = & 1- \left( \frac{1}{2} A_{-\infty} \left(1 + \erf \left( \frac{\zeta}{2} \right) \right)+ o(1) \right)  \nonumber  \\ 
&  \times \exp \left( -t  \Bigg( \bigg( \frac{ v^*(u_c) - \sqrt{v^*(u_c)^2 - 4f'(1)}}{2} \bigg)^2 -f'(1) \Bigg)  \right. \nonumber \\
&  \left. \hspace{1.32cm} - t^{\frac{1}{2}} \zeta \left( \frac{v^*(u_c) - \sqrt{v^*(u_c)^2 - 4f'(1)}}{2} \right)  \right) ,
\end{align} 
as $t \to \infty$ with $\zeta = O(1)$, and in region $\mathbf{IV_L^b}$
\beq
u(w,t) =  1 - \exp \Bigg( -t \left( \frac{v^*(u_c) - \sqrt{v^*(u_c)^2 - 4f'(1)}}{2} \right)  + \ln A_{-\infty} + o(1) \Bigg) , \eelab{4lb}
\eeq
as $t \to \infty$ with $-  \sqrt{v^*(u_c)^2 -  4 f'(1)} + O(t^{-\frac{1}{2}}) < w < O(t^{-1})^-$. 
A schematic representation of the location and thickness of the asymptotic regions as $t \to \infty$ is given in Figure \ffref{large_time}.
\begin{figure}
\centering
{\includegraphics[width=0.9\linewidth]{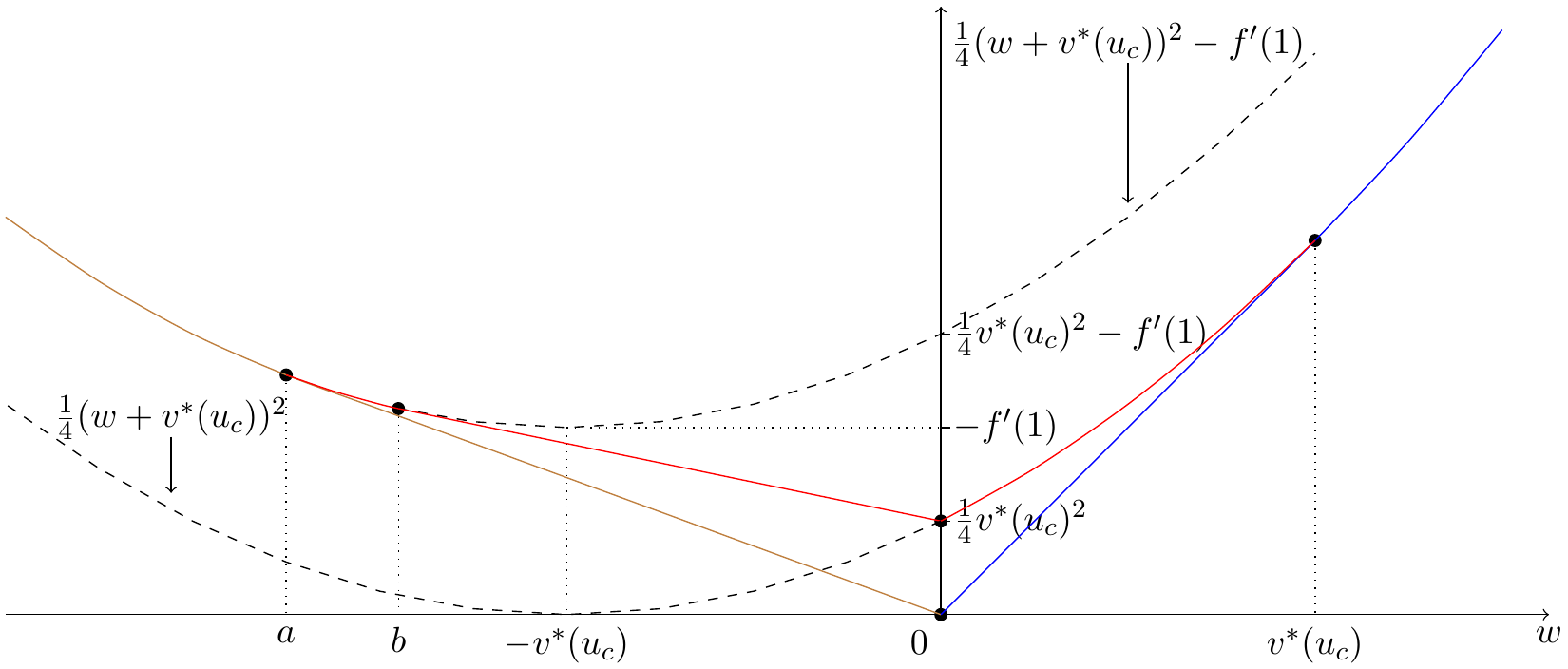}} 
\caption{Sketches of the exponent in the large-time solution to QIVP. Sketches of the leading order term $G_0(w)$ when $w < 0$  (brown), in expansions \eeref{u4la} and \eeref{4lb}, in regions $\mathbf{IV_L^a}$ and $\mathbf{IV_L^b}$, respectively; sketches of the leading order term $\bar{G}_0(w)$ when $w > 0$ (blue), in expansions \eeref{exp_u_4ra3} and \eeref{u_4rb} in regions $\mathbf{IV_R^a}$ and $\mathbf{IV_R^b}$, respectively; and sketches of the exponential corrections (red) in regions $\mathbf{IV_L^b}$ $(a < w < 0)$ and $\mathbf{IV_R^b}$ $( 0 < w < v^*(u_c))$, respectively. Here we have used the notation $a= - \sqrt{ v^*(u_c)^2 - f'(1)}$ and $b= - 2 \sqrt{- f'(1)}$.}
\fflab{region4}
\end{figure}

We next consider the structure of the expansion in region $\mathbf{T_L}$ in closer detail. Via \eeref{exp_uTL}, we observe that for $\zeta \gg 1$,
  \begin{align} \eelab{exp_uTL_match}
u(\zeta,t) \sim  1  - & \exp \left( -t  \left( \left( \frac{ v^*(u_c) - \sqrt{v^*(u_c)^2 - 4f'(1)}}{2} \right)^2 -f'(1) \right)  \right. \nonumber \\
& \left. \;\,\qquad - t^{\frac{1}{2}} \zeta \left( \frac{v^*(u_c) - \sqrt{v^*(u_c)^2 - 4f'(1)}}{2} \right)  + \ln 
\left( A_{-\infty} 
\left(1- \frac{1}{\sqrt{\pi}} \frac{1}{\zeta} e^{- \frac{\zeta^2}{4}}
\right) 
\right)\right),
\end{align} 
as $t \to \infty$, which demands that in region $\mathbf{IV_L^b}$, to continue the expansion in \eeref{4lb}, we must write
\beq
u(w,t) =  1 - A_{-\infty} \exp \Bigg[ -t \left( \frac{v^*(u_c) - \sqrt{v^*(u_c)^2 - 4f'(1)}}{2} \right) w  \Bigg] + t^{- \hat{\beta}} G(w,t) e^{-t H(w)}, \eelab{exp_u4lb}
\eeq
as $t \to \infty$ with $-  \sqrt{v^*(u_c)^2 -  4 f'(1)} + O(t^{-\frac{1}{2}}) < w < O(t^{-1})^-$ and $G(w,t)=O(1)$ as $t \to \infty$. Here $\hat{\beta}$ is a constant to be determined and  
\beq \eelab{region_4lb_1}
H(w) > \frac{1}{2} \left( v^*(u_c) - \sqrt{v^*(u_c)^2 - 4f'(1)} \right) w,
\eeq 
for all $-  \sqrt{v^*(u_c)^2 -  4 f'(1)}  < w < 0$. On substituting from expansion \eeref{exp_u4lb} with \eeref{region_4lb_1} into equation \eeref{QIVPa} we obtain
  \begin{align} 
& G \left( H_w^2 - (w+ v^*(u_c)H_w + H + f'(1) \right) + O(t^{-1} )  \nonumber \\
& \hspace{0.1cm} = O \left( t^{\gamma + \hat{\beta}} \exp \left( -t \left( \left( \frac{v^*(u_c) - \sqrt{v^*(u_c)^2 - 4f'(1)}}{2} \right) w  + \frac{1}{4}v^*(u_c)^2 - H(w) \right) \right) \right),
 \eelab{region_4lb_2}
\end{align} 
as $t \to \infty$ with $-  \sqrt{v^*(u_c)^2 -  4 f'(1)} + O(t^{-\frac{1}{2}}) < w < O(t^{-1})^-$. To obtain a non-trivial balance at leading order as $t \to \infty$ we suppose that the function $H(w)$ is such that the right-hand side of equation \eeref{region_4lb_2} is exponentially small as $t \to \infty$, and we will later verify this as consistent. Thus, at leading order, we obtain the following boundary value problem in region $\mathbf{IV_L^b}$ for $H(w)$,
\begin{subequations} \eelab{2nd_BVP_4lb}
  \begin{align}
& H_w^2 - (w+ v^*(u_c) )H_w + H =- f'(1),   \\
& 0 < H(w) - \frac{1}{2} \left( v^*(u_c) - \sqrt{v^*(u_c)^2 - 4f'(1)} \right) w <  \frac{1}{4}v^*(u_c)^2   \eelab{2nd_BVP_4lb_2},
\end{align}  
with $-  \sqrt{v^*(u_c)^2 -  4 f'(1)} < w < 0$ and which must be solved subject to 
  \begin{align}
& H(w) \to   \frac{1}{4}v^*(u_c)^2 \quad \mbox{as} \quad w \to 0^-,   \eelab{2nd_BVP_4lb_3} \\
& H(w) \sim  \frac{1}{4} (w+v^*(u_c))^2-f'(1), 
\quad \mbox{as} \quad  w \to (-  \sqrt{v^*(u_c)^2 -  4 f'(1)})^+. \eelab{2nd_BVP_4lb_5} 
\end{align}  
\end{subequations}
 Here the lower bound of inequality \eeref{2nd_BVP_4lb_2} follows from \eeref{region_4lb_1} whilst the upper bound ensures the right-hand side of equation \eeref{region_4lb_2} is exponentially small as $t \to \infty$. Condition \eeref{2nd_BVP_4lb_3} is required so that the correction term in expansion \eeref{exp_u4lb} is of the appropriate order to enable matching of \eeref{exp_u4lb} in region $\mathbf{IV_L^b}$ (as $w \to 0^-$) with expansion \eeref{exp_u5l}, \eeref{BVP_TWb},  \eeref{exp_TWb}, \eeref{region5_bal} and \eeref{region5_cor2}, in region $\mathbf{V_L}$ (as $y \to - \infty$). 
Condition \eeref{2nd_BVP_4lb_5} represents the matching condition between the expansion in region $\mathbf{IV_L^b}$ as $w \to (-  \sqrt{v^*(u_c)^2 -  4 f'(1)})^+$ (given by \eeref{exp_u4lb}) and the expansion in region $\mathbf{T_L}$ as $\zeta \to \infty$ (given by \eeref{exp_uTL_match}). Recalling that for each $u_c \in (0,1)$ then $v^*(u_c) \in (0,2)$, the boundary value problem \eeref{2nd_BVP_4lb} has the unique solution
\begin{subequations}\eelab{2nd_BVP_4lb_sol}
\beq 
H(w) = \begin{cases}  H_{L1}(w) , \quad &  -  \sqrt{v^*(u_c)^2 -  4 f'(1)} < w < -2 \sqrt{-f'(1)}  ,\\
 H_{L2}(w), \quad &  -2 \sqrt{-f'(1)}  \leq  w < 0,  \end{cases} \\
\eeq
with 
\beq
 H_{L1}(w) =  \frac{1}{4} (w+v^*(u_c))^2-f'(1)   
\quad 
\text{and}
\quad
 H_{L2}(w) = \frac{1}{4}v^*(u_c)^2 + \left( \frac{1}{2}v^*(u_c) - \sqrt{- f'(1)} \right) w,
\eeq
\end{subequations}
and where we also determine, via asymptotic matching, that $\hat{\beta} = \frac{1}{2}$ for  $- \sqrt{v^*(u_c)^2 -  4 f'(1)} + O(t^{-\frac{1}{2}}) < w < -2 \sqrt{-f'(1)} -O(t^{-\frac{1}{2}})$. A sketch of the exponents in expansions \eeref{exp_u_4ra3} and \eeref{u_4rb}, \eeref{u4la} and \eeref{4lb} in regions $\mathbf{IV_R^a},\mathbf{IV_R^b}$, $\mathbf{IV_L^a}$ and $\mathbf{IV_L^b}$, respectively, is given in Figure \ffref{region4}. We note that although $H(w)$ and $H'(w)$ are continuous for all $-  \sqrt{v^*(u_c)^2 -  4 f'(1)}  < w <0$, the second derivative $H''(w)$ is discontinuous at the point $w= -2 \sqrt{-f'(1)}$. Hence, a thin transition region about the point $w= -2 \sqrt{-f'(1)}$ is required in which the second derivative in equation \eeref{QIVPa} is   retained at leading order to smooth out the discontinuity. However, this region is passive, and for brevity will not be considered here. It remains to determine $G(w,t)$ in region $\mathbf{IV_L^b}$. To that end, since $G(w,t)=O(1)$ as $t \to \infty$ with $w=O(1)^-$, we must expand $G(w,t)$ in the form 
\beq \eelab{region4r_1}
G(w,t) = G_0(w) + t^{-\lambda} G_1(w) + o \left( t^{- \lambda} \right),
\eeq
as $t \to \infty$ with $-  \sqrt{v^*(u_c)^2 -  4 f'(1)} + O(t^{-\frac{1}{2}}) < w <O(t^{-1})$ and substitute from expansion \eeref{exp_u4lb} (with \eeref{2nd_BVP_4lb_sol} and \eeref{region4r_1}) into equation \eeref{QIVPa}. When $ -  \sqrt{v^*(u_c)^2 -  4 f'(1)} < w < -2 \sqrt{-f'(1)}$ we find $\lambda=1$ and   at leading order 
 $G_0(w)$ remains indeterminate when $ -  \sqrt{v^*(u_c)^2 -  4 f'(1)}  < w < -2 \sqrt{-f'(1)} $ and represents a further globally determined function. However, when $-2 \sqrt{-f'(1)} < w < 0$, we require that $\lambda=1$ and  at leading order we obtain 
\beq
\left( w + 2 \sqrt{- f'(1)} \right) G_0'  = - \hat{\beta} G_0 ,
\eeq
which gives, on integration, 
\beq
G_0(w) = \frac{ \left( 2 \sqrt{- f'(1)} \right)^{\hat{\beta}} A_L}{ \left( w + 2 \sqrt{-f'(1)} \right)^{\hat{\beta}}} ,
\eeq 
with $-2 \sqrt{-f'(1)} < w <  0$, where $A_L \neq 0$ is a globally determined constant. Therefore, the expansion in region $\mathbf{IV_L^b}$ is developed to,
\beq
u(w,t) =  1 - A_{-\infty} \exp \Bigg( -t \left( \frac{v^*(u_c) - \sqrt{v^*(u_c)^2 - 4f'(1)}}{2} \right) w  \Bigg) + \hat{u}(w,t) , \eelab{exp_u_4lb_full} 
\eeq
as $t \to \infty$. Here  
\beq
\hat{u}(w,t) = t^{- \beta_1} \big( G_0(w)+o(1) \big) \exp { \left( -t \left( \frac{1}{4} \left( w + v^*(u_c) \right)^2 - f'(1) \right) \right) },
\eeq
when $ -  \sqrt{v^*(u_c)^2 -  4 f'(1)} + O(t^{-\frac{1}{2}}) < w < -2 \sqrt{-f'(1)}  - O(t^{-\frac{1}{2}}) $, with
\beq
G_0(w) \sim \frac{A_{-\infty}}{\sqrt{\pi}} \left( w + \sqrt{v^*(u_c)^2 -  4 f'(1)} \right)^{-1},
\eeq
as $w \to ( -  \sqrt{v^*(u_c)^2 -  4 f'(1)} )^+$ and 
\beq
\beta_1 = \frac{1}{2},
\eeq
on matching with region $\mathbf{T_L}$. However,
 \begin{align} 
\hat{u}(w,t) = & \frac{ \left( 2 \sqrt{- f'(1)} \right)^{\beta_2} A_L}{ \left( w + 2 \sqrt{-f'(1)} \right)^{\beta_2}}  t^{- \beta_2} \big( 1 + o(1) \big) \nonumber \\
& \times \exp { \left( -t \left( \frac{1}{4}v^*(u_c)^2 + \left( \frac{1}{2}v^*(u_c) - \sqrt{- f'(1)} \right) w \right) \right) },
\end{align} 
when $-2 \sqrt{-f'(1)} + O(t^{-\frac{1}{2}}) < w <  O(t^{-1})^-$, and with $\beta_2$ undetermined at this stage. It is important to recall that the change in structure of $\hat{u}(w,t)$ across $w= -2 \sqrt{-f'(1)}$ is accommodated in a transition region when $w=-2 \sqrt{-f'(1)} \pm O(t^{- \frac{1}{2}})$. This region is passive and its details may be omitted here. 

We can now return to region $\mathbf{V_L}$. It follows from \eeref{exp_u5l} with \eeref{BVP_TWb},  \eeref{exp_TWb}, \eeref{region5_bal} and \eeref{region5_cor2}, that in region $\mathbf{V_L}$ we must develop expansion \eeref{exp_u5l} in the form 
\beq \eelab{exp_u5l_full}
u(y,t) = U_T(y) + t^{\gamma} \exp{ \left( - \frac{1}{4}  v^*(u_c)^2  t \right) } u_1(y) (1 + o(1)),
\eeq
as $t \to \infty$ with $y=O(1)^-$. On substituting from expansions \eeref{exp_large_s} and \eeref{exp_u5l_full} into equation \eeref{QIVPa}, and cancelling at leading order, we obtain 
\begin{subequations} \eelab{BVP2_5l}
 \begin{align} \eelab{BVP2_5l_1}
& u_1'' + v^*(u_c) u_1' + \left( \frac{1}{4}v^*(u_c)^2 + f'(U_T(y)) \right) u_1 = -c_3 U_T'(y), \quad y < 0 ,   \\
& u_1(0^-)=0, \eelab{BVP2_5l_2}
\end{align}  
where the condition \eeref{BVP2_5l_2} is required for the boundary condition \eeref{QIVPd} to be satisfied. It remains to match expansion \eeref{exp_u5l_full} in region $\mathbf{V_L}$ (as $y \to - \infty$) with expansion \eeref{exp_u_4lb_full} in region $\mathbf{IV_L^b}$ (as $w \to 0^-$). On applying Van Dyke's matching principle \cite{VanDyke1975}, we readily obtain this matching condition as
\beq \eelab{match_con_5l}
u_1(y) \sim A_L \exp \left(  \left( \sqrt{- f'(1)} - \frac{1}{2} v^*(u_c) \right) y \right) \quad \mbox{as} \quad y \to - \infty,
\eeq
\end{subequations}
with $\beta_2$ now determined as 
\beq
\beta_2 = - \gamma.
\eeq

On collecting \eeref{BVP2_5r} and \eeref{BVP2_5l}, in addition to the derivative continuity condition \eeref{QIVPe} at $y=0$, we obtain the following boundary value problem for $u_1(y)$,
\begin{subequations} \eelab{BVP_5_2nd}
  \begin{align} 
& u_1'' + v^*(u_c) u_1' + \frac{1}{4} v^*(u_c)^2 u_1 = -c_3 U_T'(y), \quad y > 0 ,  \eelab{BVP_5_2nd_1}   \\
& u_1'' + v^*(u_c) u_1' + \left( \frac{1}{4}v^*(u_c)^2 + f'(U_T(y)) \right) u_1 = -c_3 U_T'(y), \quad y < 0 ,  \eelab{BVP_5_2nd_2} \\
& u_1(y) \sim A_L \exp \left( \left( \sqrt{- f'(1)} - \frac{1}{2} v^*(u_c) \right) y \right) \quad \mbox{as} \quad y \to - \infty, \eelab{BVP_5_2nd_3}  \\
& u_1(0^-)=u_1(0^+) =0, \eelab{BVP_5_2nd_4}  \\
& u_1'(0^-)=u_1'(0^+), \eelab{BVP_5_2nd_5} 
\end{align} 
\end{subequations}
which must be solved subject, in addition, to the matching condition on $u_1(y)$ as $y \to \infty$ with expansion \eeref{exp_u_4rb2} in region $\mathbf{IV_R^b}$. 
We begin in $y<0$, with the inhomogeneous linear equation \eeref{BVP_5_2nd_2}. Since $U_T(y)$ satisfies the equation $U_T''(y)+v^*(u_c)U_T'(y)+f_c(U_T(y))=0$, a particular integral for \eeref{BVP_5_2nd_2} is readily deduced to be proportional to $U_T'(y)$, and so the general solution to \eeref{BVP_5_2nd_2} may be written as 
\beq \eelab{gen_sol_BVP5_2nd}
u_{1}(y) = E_0 \phi_+ (y) + E_1 \phi_- (y) -  4 \frac{c_3}{v^*(u_c)^2} U_T'(y), \quad y \leq 0,
\eeq
with $\phi_+(y), \phi_-(y): ( - \infty, 0] \to \mathbb{R}$ basis functions for the homogeneous part of equation \eeref{BVP_5_2nd_2} chosen so that 
\begin{subequations} \eelab{TW_fun}
 \begin{align}
& \phi_+(y) \sim \exp { \left(   \left( \sqrt{- f'(1)} - \frac{1}{2} v^*(u_c)       \right) y \right) }, \\
& \phi_-(y) \sim \exp { \left(  - \left( \sqrt{- f'(1)} + \frac{1}{2} v^*(u_c)       \right) y \right) },
\end{align}
\end{subequations}
as $y \to - \infty$, whilst $E_0$ and $E_1$ are arbitrary  constants to be determined. 
It follows from \eeref{exp_TWb}, \eeref{TW_fun} and an application of condition \eeref{BVP_5_2nd_3} that we must have  
\beq \eelab{E0andE1}
E_0 = A_L, \qquad E_1=0.
\eeq
Moreover, on applying condition \eeref{BVP_5_2nd_4} (where we have evaluated $U_T'(0)$ via  \eeref{exp_TWa}) we obtain 
\beq \eelab{c3}
c_3 = - \frac{A_L v^*(u_c) \phi_+(0)}{4 u_c}.
\eeq
Thus, on collecting expressions \eeref{gen_sol_BVP5_2nd}, \eeref{E0andE1} and \eeref{c3} we have  
\beq
u_{1}(y) = A_L \phi_+ (y) + \frac{A_L \phi_+(0)}{v^*(u_c) u_c } U_T'(y), \quad y < 0.
\eeq
We next consider $u_1(y)$ with $y>0$. The general solution to the inhomogeneous linear equation \eeref{BVP_5_2nd_1} (using equations \eeref{exp_TWa} and \eeref{c3}) is readily found to be 
\beq
u_1(y) = \left( E_3 + E_4 y \right) e^{- \frac{1}{2} v^*(u_c) y} - A_L \phi_+(0) e^{- v^*(u_c) y}, \quad y \geq 0,
\eeq
with arbitrary  constants $E_3$ and $E_4$ determined, via application of the coupling conditions \eeref{BVP_5_2nd_4} and \eeref{BVP_5_2nd_5}, as 
 \begin{align}
& E_3 = A_L \phi_+(0) ,  \\
& E_4 = A_L \left(  \phi_+'(0) + \phi_+(0) \left( \frac{1}{2} v^*(u_c) - \frac{f_c^+}{v^*(u_c) u_c} \right) \right) ,
\end{align}
with $A_L \neq 0$. Finally, we match the expansion in region $\mathbf{V_R}$ (as $y \to \infty$) with the expansion in region $\mathbf{IV_R^b}$ (as $w \to 0^+$). Now, when $E_4 =0$, we obtain the matching condition 
\beq
\bar{G}_0(w) \sim A_L \phi_+(0) \quad \mbox{as} \quad w \to 0^+,
\eeq
and 
\beq
\gamma = - \frac{1}{2} (= - \beta_2).
\eeq
However, when $E_4 \neq 0$, we obtain the matching condition  
\beq
\bar{G}_0(w) \sim E_4 w \quad \mbox{as} \quad w \to 0^+,
\eeq
and 
\beq
\gamma = - \frac{3}{2} (= - \beta_2).
\eeq
Also, it follows from expression \eeref{c3} (since $A_L \neq 0$) that $c_3=0$ if and only if $\phi_+(0)=0$. Therefore, we have the following cases, namely;  
\paragraph*{Case (I) $\bm{\phi_+(0) \neq 0$}.}
In this case
\beq
c_3 \neq 0, \nonumber 
\eeq
and
\beq
E_4 = 0 \; \mbox{with} \; \gamma  = - \frac{1}{2} (= - \beta_2)  \quad \mbox{or} \quad E_4 \neq 0 \; \mbox{with} \; \gamma  = - \frac{3}{2} (= - \beta_2). \nonumber  
\eeq
\paragraph*{Case (II)   $\bm{\phi_+(0) = 0}$.}
In this case $\phi_+'(0) \neq 0$ and
\beq
c_3=0, \nonumber 
\eeq
whilst $E_4 \neq 0$, and so
\beq 
\gamma  = - \frac{3}{2} (= - \beta_2).  \nonumber
\eeq
We next consider the basis function $\phi_+: (- \infty,0] \to \mathbb{R}$. For fixed $u_c \in (0,1)$ the initial value problem for $\phi_+: (- \infty,0] \to \mathbb{R}$ is given by
\begin{subequations} \eelab{BVP_base_fun_5l}
 \begin{align}
& \phi_+'' + v^*(u_c) \phi_+' + \left( \frac{1}{4}v^*(u_c)^2 + f'(U_T(y)) \right) \phi_+ = 0, \quad y < 0 , \eelab{BVP_base_fun_5l_1} \\
& \phi_+(y) \sim \exp \left( \left( \sqrt{- f'(1)} - \frac{1}{2} v^*(u_c) \right) y \right) \quad \mbox{as} \quad y \to - \infty . \eelab{BVP_base_fun_5l_2}
\end{align}
\end{subequations}
We reduce the problem \eeref{BVP_base_fun_5l} to normal form by setting $\phi_+(y)=\psi_+(y) \exp{ \left( - \frac{1}{2} v^*(u_c) y \right) } $ with $\psi_+: (- \infty,0] \to \mathbb{R}$ now satisfying the initial value problem
\begin{subequations} \eelab{BVP2_base_fun_5l}
 \begin{align}
& \psi_+'' + f'(U_T(y)) \psi_+ = 0, \quad y < 0 ,   \\
& \psi_+(y) \sim \exp \left(  \sqrt{- f'(1)}  y \right) \quad \mbox{as} \quad y \to - \infty. \eelab{BVP2_base_fun_5l_2}
\end{align}
\end{subequations}
This can now be solved numerically to find $\psi_+(0)$ and $\psi_+'(0)$ which we then use to obtain $\phi_+(0)$ and $\phi_+'(0)$, after which the occurrence of case (I) or case (II) is determined.

\bigskip

The asymptotic structure of the solution to QIVP as $t \to \infty$ is now complete with the expansions in regions $\mathbf{IV_L^a}$, $\mathbf{T_L}$, $\mathbf{IV_L^b}$, $\mathbf{V_L}$, $\mathbf{V_R}$, $\mathbf{IV_R^b}$, $\mathbf{T_R}$ and $\mathbf{IV_R^a}$ providing a uniform approximation to the solution of QIVP as $t \to \infty$. On collecting expressions \eeref{exp_large_s}, \eeref{c0}, \eeref{phi1},  \eeref{c1} and \eeref{phi3} we have obtained, in particular, that 
\beq
\dot{s}(t) \; = \; v^*(u_c) + c_3 t^{\gamma}  \exp {  \bigg( - \frac{1}{4}  v^*(u_c)^2 t \bigg) }  + o \left( t^{\gamma}  \exp {  \bigg( - \frac{1}{4}  v^*(u_c)^2 t \bigg) }  \right) \quad \mbox{as} \quad t \to \infty,  \eelab{large_time_speed}
\eeq
where the constants $c_3$ and $\gamma$ depend upon whether case (I) or case (II) is pertaining for the given KPP reaction function and the cut-off value $u_c \in (0,1)$.
Hence, via the method of matched asymptotic coordinate expansions, we have been able to obtain the correction term to the asymptotic propagation speed $v^*(u_c)$ of the developing PTW structure in the solution to QIVP as $t \to \infty$.
 In addition, with $u : \mathbb{R} \times [ 0 , \infty) \to \mathbb{R}$ being the solution to QIVP, it follows from expansions \eeref{exp_u_4ra3}, \eeref{match_con_4l1}, \eeref{u_TR}, \eeref{exp_u_4rb2}, \eeref{exp_u5r_2},  
\eeref{u4la}, \eeref{match_con_4la_full}, \eeref{exp_uTL}, \eeref{exp_u_4lb_full}, \eeref{exp_u5l_full} in regions $\mathbf{IV_L^a}$, $\mathbf{IV_L^b}$, $\mathbf{IV_R^a}$, $\mathbf{IV_R^b}$, $\mathbf{T_L}$, $\mathbf{T_R}$, $\mathbf{V_L}$ and $\mathbf{V_R}$ that,
\beq \eelab{large_time_u}
u(y,t)=U_T(y) + E(y,t),
\eeq
as $t \to \infty$ for $y \in \mathbb{R}$, with $E(y,t)$ linearly exponentially small in 
$t$ as $t \to \infty$, uniformly for $y \in \mathbb{R}$. In particular, on any closed bounded interval I,
\beq \eelab{large_time_u_corr}
E(y,t) = O \left(  t^{\gamma} e^{- \frac{1}{4} v^{*2}(u_c) t}   \right), 
\eeq
as $t \to \infty$ uniformly for $y \in $ I. A significant point to note here, is that, for KPP reaction functions satisfying 
\eeref{KPPreaction}, in the absence of cut-off, the corresponding correction terms in \eeref{large_time_speed}, \eeref{large_time_u} and \eeref{large_time_u_corr} are only algebraically small in $t$ as $t \to \infty$, being of $O(t^{-1})$ (see, for example, Leach and Needham \cite{LeachNeedham2003}).

To illustrate these results we consider a simple example of KPP reaction function $f: \mathbb{R} \to \mathbb{R}$ which satisfies \eeref{KPPreaction}, and has 
\beq
f(u) = \lambda (1-u), \quad u \geq \frac{1}{2} \left( 1 + \frac{\lambda}{(1 + \lambda)}   \right),
\eeq
with $\lambda > 0$ fixed. With the cut-off value 
\beq \eelab{example_uc}
u_c \in \left. \left[ \frac{1}{2} \left( 1 + \frac{ \lambda}{ ( 1 + \lambda) }    \right), 1 \right) \right. ,
\eeq
then, in this example, $f_c : \mathbb{R} \to \mathbb{R}$ is given by
\beq \eelab{example_reaction}
f_c(u) = \begin{cases}  0, \quad & u \in (- \infty, u_c] ,     \\
 \lambda(1-u), \quad &  u \in (u_c, \infty),  \end{cases}  
\eeq
and
\beq
f'(1) = - \lambda, \quad f_c^+ = \lambda ( 1 - u_c ).
\eeq
For this example, we can readily obtain the PTW explicitly as $U_T: \mathbb{R} \to \mathbb{R}$ given by
 \begin{align} \eelab{example_sol}
U_T(y) = \begin{cases}  1 - (1-u_c) \exp { \left(  \left( \frac{ \sqrt{ v^*(u_c)^2 + 4 \lambda } -v^*(u_c)}{2}            \right) y    \right) }, \quad & y \leq 0,     \\
 u_c e^{ - v^*(u_c) y} , \quad & y > 0,  \end{cases}  
\end{align}
with propagation speed
\beq \eelab{example_speed}
v^*(u_c) = \sqrt{  \lambda}  \frac{ ( 1 - u_c)}{\sqrt{u_c}}. 
\eeq
Now, via \eeref{BVP_base_fun_5l}, the basis function $\phi_+ : (-\infty, 0] \to \mathbb{R}$ satisfies 
\begin{subequations} \eelab{example_BVP_base_fun}
 \begin{align}
& \phi_+'' + v^*(u_c) \phi_+' + \left( \frac{1}{4}v^*(u_c)^2 - \lambda \right) \phi_+ = 0, \quad y < 0 , \eelab{example_BVP_base_fun_1}  \\
& \phi_+(y) \sim \exp \left( \left( \sqrt{ \lambda } - \frac{1}{2} v^*(u_c) \right) y \right) \quad \mbox{as} \quad y \to - \infty ,\eelab{example_BVP_base_fun_2} 
\end{align}
\end{subequations}
which has solution 
\beq \eelab{example_base_sol}
\phi_+(y) = \exp \left( \left( \sqrt{ \lambda } - \frac{1}{2} v^*(u_c) \right) y \right), \quad y \leq 0 . 
\eeq
Thus we obtain via \eeref{example_base_sol}
\beq
\phi_+(0) =1, \quad \phi_+'(0)=  \sqrt{ \lambda } - \frac{1}{2} v^*(u_c),
\eeq
and so,
\beq
E_4 = A_L \sqrt{ \lambda } \left( 1 - \sqrt{u_c} \right) \neq 0 .
\eeq
Thus, the particular reaction function \eeref{example_reaction} falls into case (I) which has 
\beq
\dot{s}(t) \;= \; v^*(u_c) + c_3 t^{-\frac{3}{2}}  \exp {  \bigg( - \frac{1}{4}  v^*(u_c)^2 t \bigg) }   + o \left( t^{- \frac{3}{2}}  \exp {  \bigg( - \frac{1}{4}  v^*(u_c)^2 t \bigg) }  \right) \quad \mbox{as} \; t \to \infty,  \eelab{example_large_time_speed}
\eeq
with $c_3 \neq 0$, and $v^*(u_c)$ given by \eeref{example_speed}. Similarly, in this example, both \eeref{large_time_u} and \eeref{large_time_u_corr} have $\gamma = - 3/2$.

 \begin{figure} [t] 
 		     \centering
 		     \begin{minipage}{0.485\textwidth}
 		  	     \includegraphics[width=\textwidth]{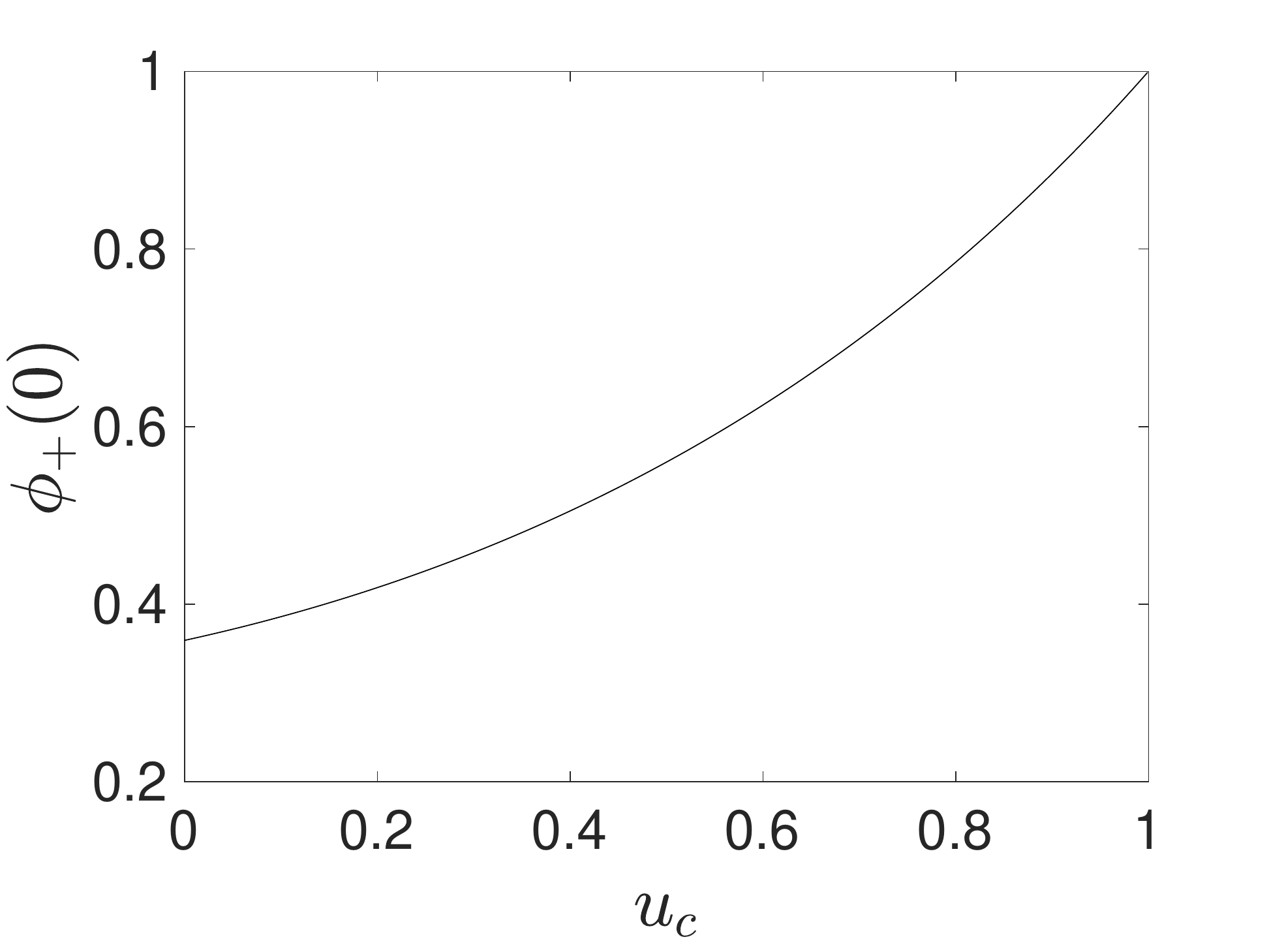}\\
 		  	      \centering  (a) 
 		     \end{minipage}
 		     \hfill
 		     \begin{minipage}{0.485\textwidth}
 		  	   \includegraphics[width=\textwidth]{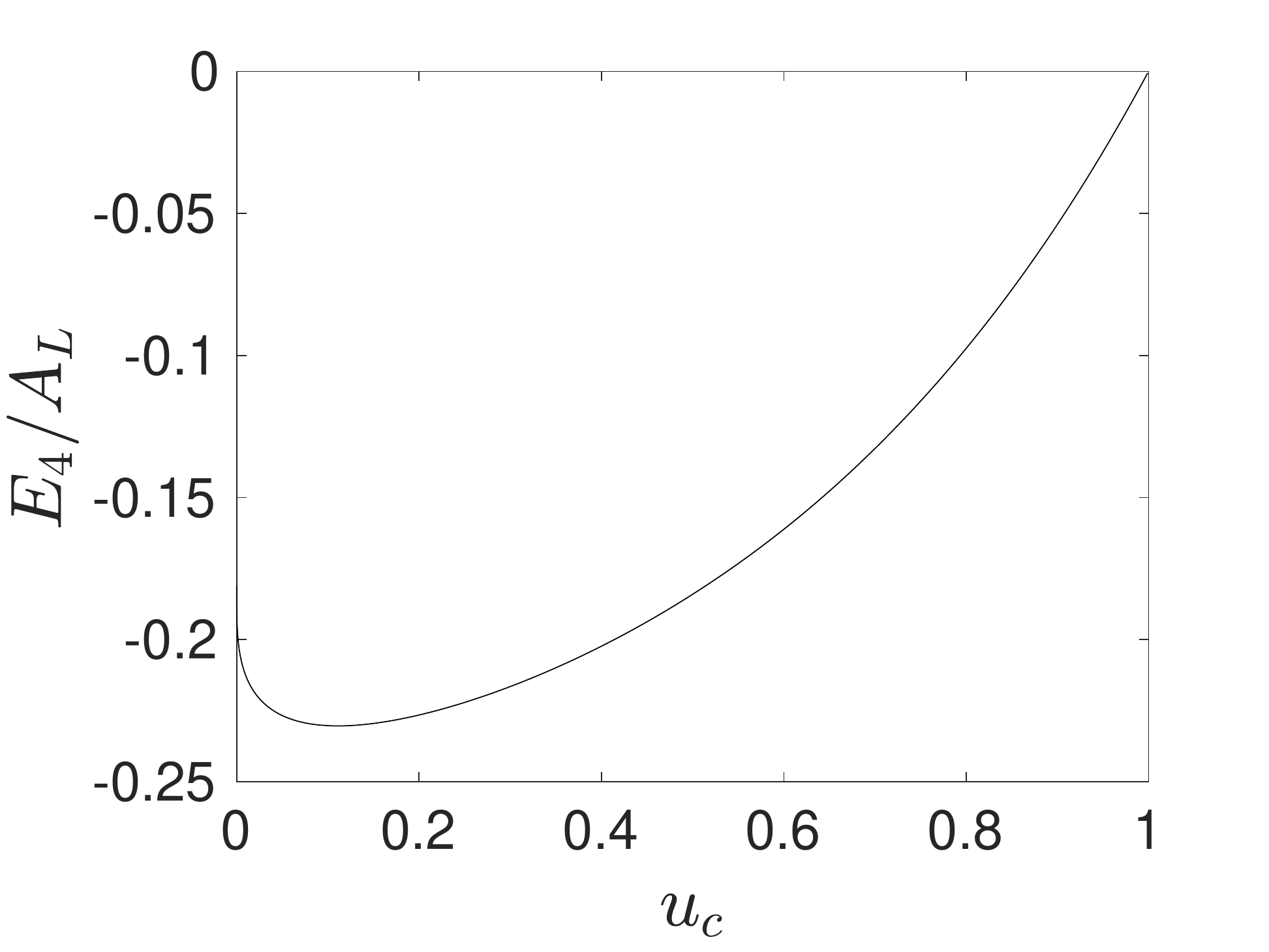}\\
 		    \centering  (b)  
			   \end{minipage}
 \caption{A graph of  (a) $\phi_+(0)$ and (b) $E_4/A_L=\phi_+'(0)+\phi_+(0)(v^*(u_c)/2-(1-u_c)/v^*(u_c))$ 
corresponding to the 
 cut-off Fisher reaction function \eeref{cutoff_Fisher}. These are obtained 
  by solving \eeref{Fisher_BVP} numerically   for a range of values of $u_c\in(0,1)$
and are used to determine the precise form of the correction to $\dot s(t)$ as $t\to\infty$, given by equation \eeref{large_time_speed}.  
 }
 \fflab{phi}
 \end{figure} 
 
 \subsection{The  case of a cut-off Fisher reaction}

To conclude this section we focus on the particular case of the cut-off Fisher reaction function \eeref{cutoff_Fisher} for fixed cut-off  $u_c \in (0,1)$.
For this example, via \eeref{BVP2_base_fun_5l}, 
  $\psi_+ : (-\infty, 0] \to \mathbb{R}$ satisfies 
\begin{subequations} \eelab{Fisher_BVP}
 \begin{align}
& \psi_+'' + \left(1 - 2 U_T(y) \right) \psi_+ = 0, \quad y < 0 , \eelab{numerics_BVP_base_fun_1}  \\
& \psi_+(y) \sim e^y \quad \mbox{as} \quad y \to - \infty. \eelab{numerics_BVP_base_fun_2} 
\end{align}
\end{subequations}
We obtain  numerical approximations of $\psi_+(0)$ and $\psi'_+(0)$ from were we deduce 
$\phi_+(0)$ and $\phi_+'(0)$.
This is readily achieved by
solving \eeref{Fisher_BVP} 
together with the nonlinear boundary value 
problem determining $U_T(y)$ (see equation (11) in Part I of this series)
numerically over an interval $y\in[-M,0]$  for $M \in\mathbb{R}^+$   using the Matlab initial value solver \texttt{ode45}, taking $v=v^*(u_c)$.
The  values of $v^*(u_c)$  and $M$ are  determined numerically as detailed in Part I of this series of papers.  
 As `initial condition'
we employ 
$(U_T,U_T',\psi_+,\psi_+')=(1-\epsilon, -\lambda_+(v^*(u_c))\epsilon, e^{-M},e^{-M})$, 
where $\epsilon=10^{-10}$ 
and prescribe an absolute and relative ODE tolerance  of $10^{-13}$.
 
Figure \ffref{phi} examines the behaviour of $\phi_+(0)$ and $E_4/A_L=\phi_+'(0)+\phi_+(0)(1/2v^*(u_c)-(1-u_c)/v^*(u_c))$ for a range of values of  $u_c$. 
It suggests that  $\phi_+(0)$ and $E_4$ are both non-zero and therefore  
the particular reaction function \eeref{cutoff_Fisher} falls into case (I) with $c_3 \neq 0$, $\gamma= - 3/2$ and where $\dot{s}(t)$ has the asymptotic expression 
\beq
\dot{s}(t) \sim  v^*(u_c) - \frac{A_L v^*(u_c) \phi_+(0)}{4 u_c} t^{-\frac{3}{2}}  \exp {  \bigg( - \frac{1}{4}  v^*(u_c)^2 t \bigg) }  \quad \mbox{as} \quad t \to \infty.  \eelab{Fisher_speed}
\eeq
 We observe that the asymptotic expression 
 \eeref{Fisher_speed} 
qualitatively  agrees with the numerical solutions for QIVP obtained for the cut-off Fisher reaction function in section \ssref{Numerical_results}:   
Figures  \ffref{sdot1} and \ffref{sdot2} suggest that 
 the correction to $\dot s(t)$ is exponentially small in $t$  as $t\to \infty$ while
  Figure \ffref{PTWu} makes clear that the exponential 
decay rate decreases with the increasing  value of $u_c$.
  However,  a quantitative test of the validity of   \eeref{Fisher_speed} 
  is challenging because 
 we do not have sufficient precision to allow the numerical solver to 
 resolve exponentially small terms in the numerical solution; as such we are unable to accurately compare  \eeref{Fisher_speed} directly with numerical solutions to estimate the global constant $A_L$.

\section{Conclusions}

In this series of papers we have considered an evolution problem for a reaction-diffusion process when the reaction function is of standard KPP-type, but experiences a cut-off in the reaction rate below the normalised cut-off concentration $u_c \in (0,1)$. We have formulated this evolution problem in terms of the moving boundary initial-boundary value problem QIVP. 
In the companion paper we considered PTW solutions $U_T(y)=u(y,t)$ to QIVP. 
In this paper  we concentrated on examining whether a PTW evolves in the large-time solution to QIVP and when this is found to be the  case, determining the rate of convergence of the solution to the PTW.
Key to this study is   $y=x-s(t)=0$ which represents the location of the moving boundary where $u=u_c$. 
We used the  method of matched asymptotic  coordinate
expansions to develop the detailed asymptotic structure of the solution to QIVP in the small-time ($t=o(1)$), intermediate-time ($t=O(1)$) and large-time ($t\to \infty$) regimes for arbitrary cut-off $u_c \in (0,1)$. 
We first determined that the asymptotic structure of $u(y,t)$ in the small-time regime
has  two   regions in $y<0$, and two  regions in $y>0$ 
and is given by  expansions \eeref{u_2l},
\eeref{par_sol_1l_full},
\eeref{par_sol_1r_full}
and
 \eeref{u_2r}.  
  The two-term asymptotic expression \eeref{exp_small_s_full} for the function $s(t)$   can be derived from the inner left and inner right regions, where $y=o(1)^-$ and $y=o(1)^+$, in addition to the leading order boundary conditions. 
This reveals that as $t \to 0^+$, $\dot{s}(t)$
has an integrable singularity which depends on the cut-off $u_c$.
Here $\dot{s}(t) \to + \infty$ when $u_c \in (0,\frac{1}{2})$, whilst, $\dot{s}(t) \to - \infty$ when $u_c \in (\frac{1}{2},1)$ with a transition case where $\dot{s}(t) \to 0$ when $u_c =\frac{1}{2}$. 
We then employed the asymptotic structure of $u(y,t)$ in the outer left and right 
 regions, where $y=O(1)^-$ and $y=O(1)^+$,  for $t=o(1)$ to determine the asymptotic structures of $u(y,t)$ when $ \lvert y \rvert \to \infty$ for  $t=O(1)$. 
The latter  is key to deriving the asymptotic structure of $u(y,t)$ 
 as 
$t\to \infty$ which consists of    two principal   regions in $y<0$ and two principal  regions in $y>0$ and given by the asymptotic expressions \eeref{exp_u_4ra3}, \eeref{match_con_4l1}, \eeref{u_TR}, \eeref{exp_u_4rb2}, \eeref{exp_u5r_2},  
\eeref{u4la}, \eeref{match_con_4la_full}, \eeref{exp_uTL}, \eeref{exp_u_4lb_full}, \eeref{exp_u5l_full}, with the asymptotic structure of $s(t)$ as $t\to \infty$ being determined simultaneously and given by the asymptotic expression \eeref{large_time_speed}.
This systematic approach 
allows to establish that the solution to QIVP 
converges to the PTW solution  
as $t\to\infty$ at a rate that   is linearly exponentially  small in $t$    
with the exact form dependent on the particular underlying KPP-type reaction function $f(u)$ and the cut-off value $u_c \in (0,1)$. 
Thus, introducing an arbitrary cut-off into the reaction   significantly modifies the    rate of convergence of the large-time solution onto the PTW (from an algebraic   to an exponential rate). Consequently, the presence of a cut-off significantly 
shortens 
  the time for the solution to QIVP to converge to the PTW. 
We anticipate that the approach developed in this paper 
will be readily adaptable to corresponding problems, when the KPP-type cut-off reaction function is replaced by a broader class of cut-off reaction functions.

\section*{Acknowledgments}
The research of A. Tisbury was supported by an EPRSC grant with reference number 1537790.

\appendix 
\section{Numerical scheme} \sslab{Numerical scheme}
We approximate $u(y,t)$ and $s(t)$ by piecewise linear functions 
$u_d(y_i,t_j)$ and  $s_d(t_j)$, 
defined on  evenly spaced  space and time grids given by 
$\{y_i =  - M + i \Delta y\}_{i=0}^{I+\mathcal{I}}$
and $\{t_j = j \Delta t\}_{j=0}^J$  
with $y_I=0$ and $t_J=T$. 
We use explicit finite differences to approximate \eeref{QIVPa} by
\beq
U_i^{j+1} - U_i^{j} = \mu \left( U_{i+1}^{j} - 2U_i^{j}+  U_{i-1}^j  \right)  + \nu\left(  S^{j+1} -  S^{j}  \right) \left(   U_{i+1}^{j} -  U_{i-1}^j  \right) 
 +  \Delta t f_c(U_i^{j}), 
\eeq
 for $i=2,\ldots,I-1,I+1,\ldots,I+\mathcal{I}-1$, 
$j=1,\ldots J$, $\mu=\Delta t/ \Delta y^2$ and $\nu=1/(2\Delta y)$, where $U_i^j=u_d(y_i,t_j)$ and $S_i^j=s_d(t_j)$ respectively approximate
  $u(y_i,t_j)$ and $s(t_j)$.
We then use \eeref{QIVPc}, \eeref{QIVPd} and \eeref{QIVPe}
to set
\beq
U^j_0=1,\quad U^j_{2I}=0, \quad U^j_I =u_c, \quad
U^j_{I+1}+ U^j_{I-1} = 2u_c, \quad\text{for $j=1,\ldots J$}.
\eeq
We solve the resulting  sparse linear algebraic system of equations 
for the unknowns 
 $U_i^j$ and $S^j$ with $i=2,\ldots,I-1,I+1,\ldots,I+\mathcal{I}-1$
 and $j=1,\ldots,J$
 in an evolutionary manner starting  from 
 \beq
 \{U^0_i\}_{i=1}^{I-1}=1,\quad \{U^0_{i}\}_{i=I}^{I+\mathcal{I}}=0,   \quad
 S^0 = 0, 
 \eeq
corresponding to the initial conditions \eeref{QIVPb} and \eeref{QIVPf}.
  We choose $\Delta y=5\times 10^{-3}$  and $\Delta t  = 0.4 \Delta y^2$ to ensure the  stability of the explicit method. 
  We 
   take $I$ and $\mathcal{I}$ sufficiently large 
   to ensure that
  any error arising from truncating the right-hand and left-hand boundary does not affect the solution in the interior. 
   In practice, we have found that choosing
   $I$ and $\mathcal{I}$ so that $e^{\lambda_+( v^*(u_c))y_0}, e^{- v^*(u_c)y_{I+\mathcal{I}}} \lesssim  5\times 10^{-5}$
   (corresponding to the asymptotic behaviour of the PTW as described by equation \eeref{exp_TW})
	provides   reasonable
	 accuracy.
 Comparison with results obtained for a spatial resolution of $\Delta y= 10^{-3}$ resulted in a less than $0.5\%$ difference in $u_d(y_i,t_j)$ and $s_d(t_j)$.
%

\bibliographystyle{abbrv}
\providecommand{\noopsort}[1]{}\providecommand{\singleletter}[1]{#1}%

\end{document}